\documentclass[10pt]{elsarticle}
\usepackage{graphicx}
\usepackage{epstopdf}
\ifpdf
  \DeclareGraphicsExtensions{.eps,.pdf,.png,.jpg}
\else
  \DeclareGraphicsExtensions{.eps}
\fi

\usepackage[margin=0.8in]{geometry}
\usepackage[utf8]{inputenc}

\usepackage{xcolor}
\definecolor{darkolivegreen}{rgb}{0.33, 0.42, 0.18}
\definecolor{celestialblue}{rgb}{0.29, 0.59, 0.82}

\usepackage[colorlinks,
            linkcolor=darkolivegreen,
            citecolor=darkolivegreen,
            urlcolor=darkolivegreen]{hyperref}
\usepackage{url}
\makeatletter
\g@addto@macro{\UrlBreaks}{\UrlOrds}
\makeatother

\hypersetup{pdftitle={CGTJS}, pdfauthor={Hern\'an Mella Lobos},}

\usepackage{fullpage}
\usepackage{amsmath, amsthm, amssymb}
\usepackage{setspace}
\usepackage{lineno}

\usepackage{booktabs}
\usepackage{subfig}
\usepackage{adjustbox}
\usepackage{float}

\title{A mixed and hybrid PML formulation for the 2D three-field dynamic poroelastic equations} 

\author[UCVLabel]{Hernán Mella\corref{cor1}}
\author[UCLabel]{Esteban Sáez}
\author[USMLabel]{Joaquín Mura}

\affiliation[UCVLabel]{organization={School of Electrical Engineering, Pontificia Universidad Católica de Valparaíso, Valparaíso},country={Chile}}
\affiliation[UCLabel]{organization={Department of Structural and Geotechnical Engineering, Pontificia Universidad Católica de Chile, Santiago},country={Chile}}
\affiliation[USMLabel]{organization={Department of Mechanical Engineering, Universidad Técnica Federico Santa María, Santiago}, country={Chile}}

\cortext[cor1]{hernan.mella@pucv.cl. Av. Brasil 2147, Valparaíso 2340000, Región de Valparaíso, Chile.}

\def\[{\left[}
\def\]{\right]}
\def\<{\langle}
\def\>{\rangle}
\def\({\left(}
\def\){\right)}
\def\[{\left [}
\def\]{\right]}
\def\({\left(}
\def\){\right)}

\newcommand{\e}[1]{\times 10^{#1}}

\newcommand{\uu}{\boldsymbol{u}}
\newcommand{\ww}{\boldsymbol{w}}
\newcommand{\St}{\mathbf{S}}
\newcommand{\Et}{\mathbf{E}}
\newcommand{\et}{\mathbf{e}}

\newcommand{\lp}{\lambda_p}

\newcommand{\tuu}{\tilde{\boldsymbol{u}}}
\newcommand{\tww}{\tilde{\boldsymbol{w}}}
\newcommand{\tSt}{\tilde{\mathbf{S}}}

\newcommand{\tlp}{\tilde{\lambda}_p}
\newcommand{\tp}{\tilde{p}}
\newcommand{\tpi}{\tilde{\pi}}

\newcommand{\scalarinnerRD}[2]{\int_{\OmegaRD}#1~ #2~d\Omega}
\newcommand{\scalarinnerPML}[2]{\int_{\OmegaPML}#1~ #2~d\Omega}
\newcommand{\innerRD}[2]{\int_{\OmegaRD}#1\cdot #2~d\Omega}
\newcommand{\innerPML}[2]{\int_{\OmegaPML}#1\cdot  #2~d\Omega}
\newcommand{\tensorinnerRD}[2]{\int_{\OmegaRD}#1 : #2~d\Omega}
\newcommand{\tensorinnerPML}[2]{\int_{\OmegaPML}#1 : #2~d\Omega}

\newcommand{\scalarinnerI}[2]{\int_{\GammaI}#1~ #2~d\Gamma}
\newcommand{\innerNRD}[2]{\int_{\GammaNRD}#1\cdot #2 d\Gamma}

\newcommand{\innerg}[2]{\int_{\Gamma_g}#1\cdot #2 d\Gamma}

\newcommand{\RR}{\mathbb{R}} 

\newcommand{\OmegaRD}{\Omega^{\text{RD}}}
\newcommand{\OmegaPML}{\Omega^{\text{PML}}}

\newcommand{\GammaI}{\Gamma_I}
\newcommand{\GammaN}{\Gamma_N}
\newcommand{\GammaNRD}{\Gamma_N^{\text{RD}}}
\newcommand{\GammaNPML}{\Gamma_N^{\text{PML}}}

\newcommand{\GammaDPML}{\Gamma_D^{\text{PML}}}

\newcommand{\tsigma}{\boldsymbol{\sigma}}
\newcommand{\tsigmaPML}{\boldsymbol{\sigma}^{\text{PML}}}
\newcommand{\tidentity}{\textbf{I}}
\newcommand{\Lambdae}{\boldsymbol{\Lambda}_e}
\newcommand{\Lambdap}{\boldsymbol{\Lambda}_p}
\newcommand{\Lambdaee}{\tilde{\boldsymbol{\Lambda}}_e}
\newcommand{\Lambdapp}{\tilde{\boldsymbol{\Lambda}}_p}

\newcommand{\operator}{\mathcal{J}}

\begin{document}

\begin{abstract}

Simulation of wave propagation in poroelastic half-spaces presents a common challenge in fields like geomechanics and biomechanics, requiring Absorbing Boundary Conditions (ABCs) at the semi-infinite space boundaries. Perfectly Matched Layers (PML) are a popular choice due to their excellent wave absorption properties. However, PML implementation can lead to problems with unknown stresses or strains, time convolutions, or PDE systems with Auxiliary Differential Equations (ADEs), which increases computational complexity and resource consumption.

This article presents two new PML formulations for arbitrary poroelastic domains. The first formulation is a fully-mixed form that employs time-history variables instead of ADEs, reducing the number of unknowns and mathematical operations. The second formulation is a hybrid form that restricts the fully-mixed formulation to the PML domain, resulting in smaller matrices for the solver while preserving governing equations in the interior domain. The fully-mixed formulation introduces three scalar variables over the whole domain, whereas the hybrid form confines them to the PML domain.

The proposed formulations were tested in three numerical experiments in geophysics using realistic parameters for soft sites with free surfaces. The results were compared with numerical solutions from extended domains and simpler ABCs, such as paraxial approximation, demonstrating the accuracy, efficiency, and precision of the proposed methods. The article also discusses the applicability of these methods to complex media and their extension to the Multiaxial PML formulation.

The codes for the simulations are available for download from \url{https://github.com/hmella/POROUS-HYBRID-PML}.
\end{abstract}


\maketitle

\newpage
\section{Introduction}
Fluid-saturated porous media are a common occurrence in nature. For instance, soils and rocks are often saturated with water in practical cases, while living tissues are saturated with blood and air. In both cases, if the solid skeleton's displacements and strains are relatively small, linear elasticity provides an accurate representation of the underlying dynamics. Additionally, when loads are applied quickly and inertial forces play a significant role, a proper modeling strategy for wave propagation in poroelastic media is necessary. As noted by Zienkiewicz et al. \cite{Zienkiewicz1980}, Biot's poroelastic theory can be employed to describe wave propagation in poroelastic media, such as in problems of traffic-induced vibrations or geophysical applications involving seismic wave propagation. The main challenge in these types of problems is properly handling outgoing waves. In the directions where outgoing waves travel, finite energy considerations lead to the so-called ``radiation conditions" towards infinity, which are used by Integral Equations or Boundary-Element methods to determine Green kernels and solve the problem rigorously. However, these conditions are often difficult to calculate and are limited to homogeneous and isotropic material properties at infinity \cite{bem2018}. An alternative solution is to use a foam-like subdomain to confine the region of interest, creating virtual windows in space to focus computational efforts on a specific area of the problem. The subdomain must effectively absorb the outgoing waves from the virtual window.

Numerous numerical methods have been proposed as Absorbing Boundary Conditions (ABCs). Local ABCs are often used for dry elastic problems or single-phase media due to their ease of implementation and local character in both time and space \cite{kausel_local_1988}. However, fluid-saturated porous media or two-phase media presents a different challenge due to an interaction between the solid skeleton and fluid flow, which depends on the loading rate. According to Biot's theory, high-frequency loading generates two dilatational waves and one shear wave. When the porous media has low permeability and the loading is within the low frequency range, the fluid's relative motion with respect to the soil is negligible, and viscous coupling dampens out the second dilatational wave \cite{biot_theory_1956,DudleyWard2017}. In this case, the fluid-saturated porous media behaves like a single-phase medium, where only one dilatational and one shear wave propagate.

In recent decades, the Perfectly Matched Layer (PML) has gained popularity as an ABC due to its excellent energy-absorbing properties. PML was first developed by Berenger \cite{Berenger1994} in the context of electromagnetism. Although the initial development was for Maxwell's equations, its use as an ABC for acoustic \cite{qi:1998}, elastic \cite{chew_perfectly_1996}, and poroelastic \cite{zeng:2001} domains was later extended. Since then, the technique has been widely used to simulate the propagation of elastic \cite{kuzuoglu_frequency_1996, wang_finite-difference_2003, basu_perfectly_2004, correia_development_2005, kucukcoban:2011, kucukcoban:2013, Francois2021a, Zhou2016} and poroelastic waves \cite{song_application_2005,martin_unsplit_2008,He2019,He2019a} and new and novel formulations have been introduced. These forms can be divided into split-field and unsplit-field approaches, both of which have drawbacks. Split-field formulations often result in mixed problems where stresses or strains are unknowns, increasing the computational cost of solving the problem \cite{drossaert_complex_2007,zeng:2001,zeng_staggered-grid_2001,ezziani_modelisation_2005}. Unsplit-field formulations typically require the estimation of convolutions or solving Auxiliary Differential Equations (ADEs) \cite{martin_unsplit_2008, He2019}, which can also be expensive due to the increased number of mathematical operations or the introduction of auxiliary variables. Additionally, little attention has been paid to simulating poroelastic waves in arbitrary domains with realistic subsoil properties.

In this article, we propose two new formulations of the Perfectly Matched Layer (PML) method for the second order three-field Biot's equations to address the previously mentioned limitations. Our fully-mixed and hybrid formulations maintain the second-order in time structure of the original equations, which makes them compatible with most time integration schemes. Furthermore, both methods only introduce three additional scalar variables, which is at least 50\% less computationally expensive than previous developments \cite{He2019}. The hybrid formulation modifies Biot's equations only in the PML region, resulting in significant computational cost savings. Our proposed methods perform well under challenging conditions, such as free-surface wave propagation, transitions between water and air-filled soft media, and complex geometries, making them suitable for simulating realistic media.

\section{Poro-elastodynamic equations}


The three-field model proposed by Biot \cite{Biot1962,Zienkiewicz1980,zienkiewicz_dynamic_1984} considers a domain $\Omega \subseteq \RR^2$ where the solid displacement $\uu(\boldsymbol{x},t)$, the displacement of the fluid phase relative to the solid $\ww(\boldsymbol{x},t)$, and the pore pressure in the fluid $p(\boldsymbol{x},t)$ interacts at any position $\boldsymbol{x}$ and time $t$ such that $(\boldsymbol{x},t) \in \Omega\times T$ with $T= \left( 0,\infty \right)$, according to:
\begin{subequations}
    \label{eq: Biot's equations}
    \begin{alignat}{3}
            \rho \, \Ddot{\uu} + \rho_f \, \Ddot{\ww} &= \nabla\cdot\tsigma && \quad\text{in }\Omega\times T
            \label{eq: Biot's equations (a)}
            \\
            \rho_f \,\Ddot{\uu} + \rho_w \, \Ddot{\ww} + \frac{\eta}{\kappa}\Dot{\ww} &= -\nabla p  && \quad\text{in }\Omega\times T
            \label{eq: Biot's equations (d)}
            \\
            -\Dot{p} &= \nabla\cdot\left\{M(\alpha\Dot{\uu} + \Dot{\ww})\right\}  && \quad\text{in }\Omega\times T
            \label{eq: Biot's equations (e)}
    \end{alignat}    
\end{subequations}
where the effective density $\rho=\rho_s(1-\phi) + \rho_f\phi \;$ depends on the solid $(\rho_s)$ and fluid densities $(\rho_f)$, as well as the porosity $\phi$. The fluid density $\rho_w=\tau \rho_f/\phi$ depends on the tortuosity $\tau$, the dynamic viscosity of the fluid $\eta$, and the saturated permeability of the porous media $\kappa$. The parameter $\alpha$ represents the Biot-Willis coefficient and $M$ the fluid-solid coupling bulk modulus, defined as
\begin{equation}
    \alpha = 1 - \frac{K_b}{K_s},\quad M = \left(\frac{\phi}{K_f}+\frac{\alpha-\phi}{K_s}\right)^{-1}
\end{equation}
where $K_b$, $K_s$, $K_f$ are the bulk moduli of the dry porous skeleton, solid, and fluid, respectively.

The constitutive law for the poroelastic media in \eqref{eq: Biot's equations} is given by
\begin{subequations}
    \begin{alignat}{2}
        \tsigma(\uu,p) &= C \et(\uu) - \alpha p \tidentity \\
        \et(\uu) &= \frac{1}{2}\left\{ \nabla\uu+(\nabla\uu)^T \right\}
    \end{alignat}    
\end{subequations}
where $C$ is the fourth-order elastic tensor with components $C_{ijkl}=\lambda_b\delta_{ij}\delta_{kl}+\mu_b(\delta_{ik}\delta_{jl} + \delta_{il}\delta_{jk})$ (with $\delta_{ij}$ the Kronecker delta), $\et$ the linear strain tensor, and $\tidentity$ the identity tensor.

The domain $\Omega$ is intended to be an unbounded semi-space, with the top boundary consisting of two parts: $\Gamma=\GammaN\cup\Gamma_g$. The free surface, denoted by $\GammaN$, is the region where tractions must vanish \cite{Morency2008}. On the other hand, $\Gamma_g$ is the part of the top boundary where an external load $\boldsymbol{g}$ is being applied. Therefore, the boundary conditions of the problem are:
\begin{subequations}
    \label{eq: boundary conditions}
    \begin{alignat}{2}
        \tsigma(\uu,p)\cdot\boldsymbol{n} &= \boldsymbol{g} \quad \text{in } \Gamma_g\times T \\
        \tsigma(\uu,p)\cdot\boldsymbol{n} &= \boldsymbol{0} \quad \text{in } \Gamma_N\times T \\
        p &= 0 \quad \text{in } \GammaN\times T
    \end{alignat}
\end{subequations}
The variables $\uu$, $\ww$, and $p$ must vanish as the distance from the source tends to infinity.

\section{Derivation of PML formulas}

\begin{figure}[h]
    \centering
    \subfloat[]{\includegraphics[scale=0.4]{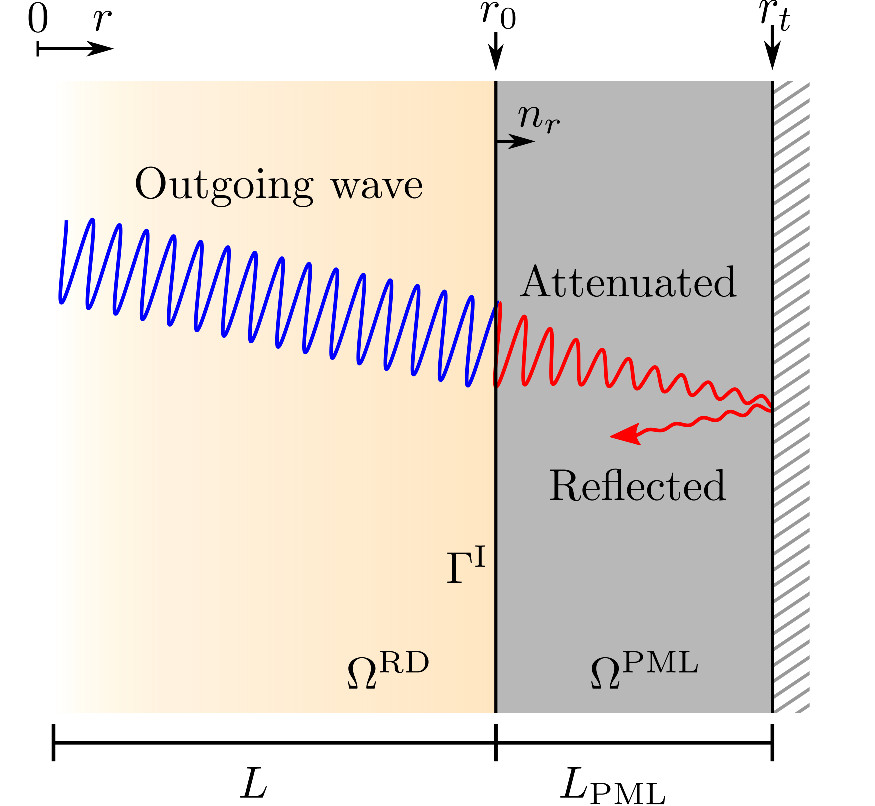} \label{fig: PML truncation}} \hfil
    \subfloat[]{\includegraphics[scale=0.55]{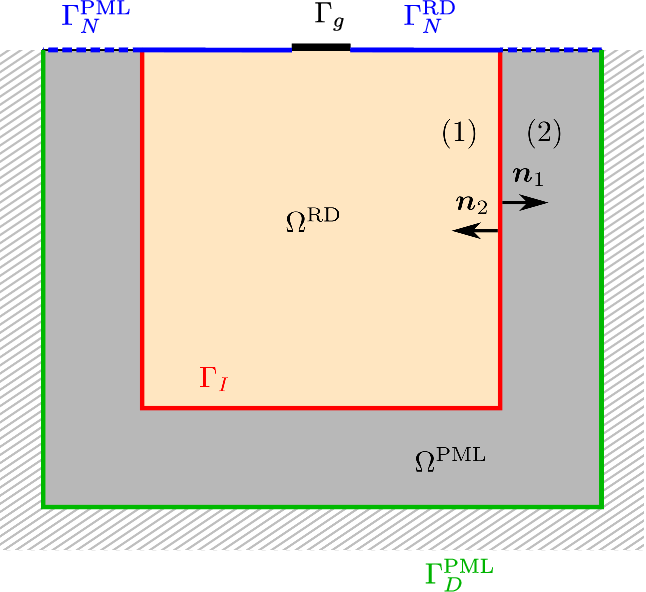} \label{fig: truncated semispace PML}}    
    \caption{(a) Semi-infinite domain $\Omega$ and (b) truncated regular and PML domains $\OmegaRD$ and $\OmegaPML$ respectively.} 
    \label{fig: PML reference domain}
\end{figure}

The region of interest where boundary-reflected waves are not desired (namely Regular Domain) $\OmegaRD$ is surrounded and truncated by a thin Perfectly Matched Layer $\OmegaPML$ such that $\Omega=\OmegaRD\cup\OmegaPML$ (see Figure \ref{fig: truncated semispace PML}). Thus, $\OmegaPML$ is defined as an extension of $\OmegaRD$ where the outgoing waves are attenuated by a complex coordinate stretching (see Figure \ref{fig: PML reference domain}).

\subsection{Complex-coordinate stretching}  

A complex-coordinate stretching applied to Biot's equations \eqref{eq: Biot's equations} leads to a modified set of equations within $\Omega$. Thus, the complex-coordinate stretching follows:
\begin{equation}
    \label{eq: stretching function}
    r \longmapsto \int_0^{r}\varepsilon_r(r',s)~dr'
\end{equation}
where $r$ denotes the spatial coordinate being transformed (namely $x$ or $y$ in the two-dimensional case), $s$ the dual variable of the time in the Laplace domain, and $\varepsilon_r$ a complex-coordinate stretching function along the coordinate $r$. The coordinate transformation given in \eqref{eq: stretching function} implies:
\begin{equation}
    \label{eq: PML derivative}
    \frac{\partial}{\partial r} 
    \longmapsto
    \frac{1}{\varepsilon_r(r)}\frac{\partial}{\partial r}
\end{equation}
which is the fundamental relation used to transform the governing equations.

The function $\varepsilon_r$ in \eqref{eq: PML function} can adopt diverse forms to achieve different absorption properties. For instance, Kuzuoglu and Mitra \cite{kuzuoglu_frequency_1996} introduced the Convolutional Frequency Shifted (CFS) PML method by defining a frequency-dependent stretching function to improve the absorption efficiency at different frequencies and also to improve the time stability. Later, Correia and Jin \cite{correia_development_2005} introduced the higher-order PML with a better absorption rate than CFS-PML. However, both approaches make the real and imaginary parts of $\varepsilon_r$ frequency-dependent, leading to convolution terms in the PML formulation. Later, Meza-Fajardo and Papageorgiou developed the Multiaxial PML (M-PML) method \cite{meza-fajardo_nonconvolutional_2008} and introduced multi-directional attenuation functions with better time stability properties. Recently, François et al. \cite{Francois2021a} proposed a non-convolutional version of the CFS-PML method in elastodynamics by introducing auxiliary variables. A standard selection for $\varepsilon_r$ would be:
\begin{equation}
    \label{eq: PML function}
    \varepsilon_r(r,s) = \alpha_r(r)+\frac{\beta_r(r)}{s},
\end{equation}
with $\alpha_r$ and $\beta_r$ denoting real-valued scaling and attenuation functions, respectively. 

The real part of $\varepsilon_r$ scales the spatial coordinate $r$ while the imaginary part is responsible for the amplitude decay of the waves entering the PML region (see Figure \ref{fig: PML reference domain}). To avoid modifying the propagating waves inside $\OmegaRD$ and ensure the wave attenuation inside $\OmegaPML$, the following conditions must be fulfilled: $\alpha_r(r)$ and $\beta_r(r)$ are constant inside $\OmegaRD$ and take values of $1$ and $0$ respectively, and both functions increases monotonically with $r$ within $\OmegaPML$. Several ways of defining $\alpha_r$ and $\beta_r$ have been proposed in the literature in different contexts \cite{Berenger1994,collino_application_2001,kucukcoban:2011,He2019}, but in this investigation, polynomial profiles were chosen according to:
\begin{subequations}
    \label{eq: PML polynomial profiles}
    \begin{alignat}{2}
        \alpha_r(r) &= \left\{\begin{array}{cl}
            1 & \text{if } 0\leq r \leq r_0 \\
            1 + \alpha_0\left\{\frac{(r-r_0)n_r}{L_{\text{PML}}}\right\}^m & \text{if } r_0 \leq r \leq r_t
        \end{array} \right. \\
        \beta_r(r) &= \left\{\begin{array}{cc}
            0 & \text{if } 0\leq r \leq r_0 \\
            \beta_0\left\{\frac{(r-r_0)n_r}{L_{\text{PML}}}\right\}^m & \text{if } r_0 \leq r \leq r_t
        \end{array} \right.
    \end{alignat}
\end{subequations}
where $n_r$ denotes the $r$-th component of the outward normal to the interface between $\OmegaRD$ and $\OmegaPML$, $L_{\text{PML}}$ the width of the PML layer, and $m$ the order of the attenuation profiles, $r_0$ and $r_t$ represent the start and end of the absorbing layer (see Figure \ref{fig: PML reference domain}). The constants $\alpha_0$ and $\beta_0$ define the absorption rate in $\OmegaPML$, and can be chosen as \cite{kucukcoban:2011}:
\begin{equation}
    \label{eq: PML scaling constants}
    \alpha_0  = \frac{(m+1)b}{2L_{\text{PML}}}\log\left(\frac{1}{|R|}\right), \qquad 
    \beta_0  = \frac{(m+1)c}{2L_{\text{PML}}}\log\left(\frac{1}{|R|}\right)
\end{equation}
where $b$ denotes a characteristic length (e.g., the width of the distributed load applied in the Neumann boundary or the cell size of the finite element mesh), $c_p$ the propagation velocity of the fastest wave, and $R$ a reflection coefficient.

\subsection{PML derivation in the Laplace domain}

The complex-coordinate stretching is enforced by introducing the relation \eqref{eq: PML derivative} into the Laplace-transformed motion equations. Applying the Laplace transform to Biot's equations leads:
\begin{subequations}
    \label{eq: Laplace Biot's equations}
    \begin{alignat}{3}     
        s^2\rho\hat{\uu} + s^2\rho_f\hat{\ww} &= \nabla\cdot\hat{\tsigma} & \text{(linear momentum conservation)} \label{eq: Laplace Biot's equations (a)}
        \\
        s^2\rho_f\hat{\uu} + s^2\rho_w\hat{\ww} + s\frac{\eta}{\kappa}\hat{\ww} &= -\nabla \hat{p}  
        & \text{(Darcy law)} \label{eq: Laplace Biot's equations (d)}
        \\
        -s\hat{p} &= \nabla\cdot\left\{M(\alpha s\hat{\uu} + s\hat{\ww})\right\}
        \label{eq: Laplace Biot's equations (e)}
        \\
        \hat{\tsigma} &= C\hat{\et} - \alpha \hat{p} \tidentity 
        & \text{(constitutive relations)} \label{eq: Laplace Biot's equations (b)}
        \\
        \hat{\et} &= \frac{1}{2}\left\{ \nabla\hat{\uu} + (\nabla\hat{\uu})^T \right\}
        \label{eq: Laplace Biot's equations (c)}
    \end{alignat}
\end{subequations}
where $\hat{f}$ denotes the variable $f$ in the Laplace domain.

\subsubsection{Linear momentum conservation}

 Replacing \eqref{eq: PML derivative} into \eqref{eq: Laplace Biot's equations (a)} for the plane strain case gives:
\begin{subequations}
    \begin{alignat}{2}
        s^2\rho\hat{u}_{x} + s^2\rho_f\hat{w}_{x} &= \frac{1}{\varepsilon_{x}}\frac{\partial \sigma_{xx}}{\partial x} + \frac{1}{\varepsilon_{y}}\frac{\partial \sigma_{xy}}{\partial y} - \frac{1}{\varepsilon_{x}}\frac{\partial \alpha \hat{p}}{\partial x} \\
        s^2\rho\hat{u}_{y} + s^2\rho_f\hat{w}_{y} &= \frac{1}{\varepsilon_{x}}\frac{\partial \sigma_{yx}}{\partial x} + \frac{1}{\varepsilon_{y}}\frac{\partial \sigma_{yy}}{\partial y} - \frac{1}{\varepsilon_{y}}\frac{\partial \alpha \hat{p}}{\partial y}
    \end{alignat}
\end{subequations}
which after multiplying both equations by $\varepsilon_x\varepsilon_y$, introducing the variables $a=\alpha_x\alpha_y$, $b=\alpha_x\beta_y + \alpha_y\beta_x$, $c=\beta_x\beta_y$, and rearranging terms can be rewritten as:
\begin{equation}
    \label{eq: Laplace PML 1}
        (s^2 a + sb + c)(\rho\hat{\uu} + \rho_f\hat{\ww}) 
        = \nabla\cdot\left\{\tsigma(\uu,p)
        \left( \Lambdae + \frac{1}{s} \Lambdap \right) 
        \right\}
\end{equation}
where the tensors $\Lambdae$ and $\Lambdap$ are defined as
\begin{equation}
    \label{eq: Lambda tensors}
    \left[\begin{array}{cc}
        \alpha_y & 0  \\
        0 & \alpha_x 
    \end{array}\right] + \frac{1}{s} \left[\begin{array}{cc}
        \beta_y & 0  \\
        0 & \beta_x 
    \end{array}\right] = \Lambdae + \frac{1}{s}\Lambdap
\end{equation}

\subsubsection{Darcy equation}
Proceeding similarly with \eqref{eq: Laplace Biot's equations (d)} results in:
\begin{subequations}
    \begin{alignat}{2}
        s^2\rho_f\hat{u}_x + s^2\rho_w\hat{w}_x + s\frac{\eta}{\kappa}\hat{w}_x &= -\frac{1}{\varepsilon_x}\frac{\partial \hat{p}}{\partial x} \\
        s^2\rho_f\hat{u}_y + s^2\rho_w\hat{w}_y + s\frac{\eta}{\kappa}\hat{w}_y &= -\frac{1}{\varepsilon_y}\frac{\partial \hat{p}}{\partial y}
    \end{alignat}
\end{subequations}
which after multiplying the first equation by $\varepsilon_x$ and the second by $\varepsilon_y$ lead to the following modified equation:
\begin{equation}
    \label{eq: Laplace PML 2}
    (\Lambdaee s + \Lambdapp)\left(s\rho_f\hat{\uu} + s\rho_w\hat{\ww} + \frac{\eta}{\kappa}\hat{\ww}\right) = -\nabla\hat{p},
\end{equation}
where the tensor $\Lambdaee$ and $\Lambdapp$ are defined as:
\begin{equation}
    \left[\begin{array}{cc}
        \alpha_x & 0  \\
        0 & \alpha_y 
    \end{array}\right] + \frac{1}{s} \left[\begin{array}{cc}
        \beta_x & 0  \\
        0 & \beta_y 
    \end{array}\right] = \Lambdaee + \frac{1}{s}\Lambdapp,
\end{equation}

The tensors $\Lambdaee$ and $\Lambdapp$ are similar to $\Lambdae$ and $\Lambdap$ (respectively) but have their diagonal reversed.

\subsubsection{Constitutive relations}

Finally, multiplying the equations \eqref{eq: Laplace Biot's equations (e)} by $\varepsilon_x\varepsilon_y$ and \eqref{eq: Laplace Biot's equations (c)} by $s\varepsilon_x\varepsilon_y$, and rearranging terms, we obtain the following PML equations in the Laplace domain:
\begin{subequations}
    \label{eq: Laplace PML}
    \begin{alignat}{2}
        \label{eq: Laplace PML 3a}
        s a\hat{\et} + b\hat{\et} + \frac{1}{s}c\hat{\et} &= \frac{1}{2}s\left\{ \nabla\hat{\uu}\Lambdae + (\nabla\hat{\uu}\Lambdae)^T \right\} + \frac{1}{2}\left\{\nabla\hat{\uu}\Lambdap + (\nabla\hat{\uu}\Lambdap)^T \right\} \\
        \label{eq: Laplace PML 3b}
        -\left(sa\hat{p} + b\hat{p} + \frac{1}{s}c\hat{p}\right) &= \nabla\cdot\left\{M(\Lambdae 
        s + \Lambdap)(\alpha\hat{\uu} + \hat{\ww})\right\}   
    \end{alignat}
\end{subequations}


\subsection{Fully-mixed time-domain formulation of the PML equations}

Applying the inverse Laplace transform to the Equations \eqref{eq: Laplace PML 1}, \eqref{eq: Laplace PML 2}, \eqref{eq: Laplace PML 3a}, and \eqref{eq: Laplace PML 3b} we obtain the time domain PML formulation of the Biot's equations:
\begin{subequations}
    \label{eq: time domain PML formulation}
    \begin{alignat}{2}
         \rho(a\Ddot{{\uu}} + b\Dot{\uu} + c\uu) + \rho_f(a\Ddot{{\ww}} + b\Dot{\ww} + c\ww) &= \nabla\cdot\left( \sigma \Lambdae + \int_0^t \sigma d\tau \Lambdap \right) \\
        a\Dot{\et} + b\et + c\int_0^t\et~d\tau &= \frac{1}{2}\left\{\nabla\Dot{\uu}\Lambdae + (\nabla\Dot{\uu} \Lambdae)^T + \nabla\uu\Lambdap + (\nabla\uu \Lambdap)^T \right\} \\
        -\nabla p &= \left(\Lambdaee\frac{\partial}{\partial t} + \Lambdapp\right)\left(\rho_f\Dot{\uu} + \rho_w\Dot{\ww} + \frac{\eta}{\kappa}\ww\right) \\
        -\left(a\Dot{p} + b p + c\int_0^t p~d\tau\right) &= \nabla\cdot\left\{M\left(\Lambdae\frac{\partial}{\partial t} + \Lambdap\right)(\alpha\uu + \ww)\right\}
    \end{alignat}
\end{subequations}
To avoid using auxiliary differential equations or the discrete evaluation of time integrals, we introduce the auxiliary memory variables $\St(\boldsymbol{x},t)$, $\Et(\boldsymbol{x},t)$, and $\pi(\boldsymbol{x},t)$ for the stress, strain, and pressure (respectively), defined as:
\begin{equation}
    \St(\boldsymbol{x},t) = \int_0^t C\et(\boldsymbol{x},\tau)~d\tau, \quad \Et(\boldsymbol{x},t) = \int_0^t\et(\boldsymbol{x},\tau)~d\tau, \quad \pi(\boldsymbol{x},t) = \int_0^t p(\boldsymbol{x},\tau)~d\tau,
\end{equation}
Consequently
\begin{subequations}
    \label{eq: mixed variables}
    \begin{alignat}{3}
        \Dot{\St}(\boldsymbol{x},t) &= C\et(\boldsymbol{x},t),\quad &&\Ddot{\St}(\boldsymbol{x},t) &&= C\Dot{\et}(\boldsymbol{x},t), \\
        \Dot{\Et}(\boldsymbol{x},t) &= \et(\boldsymbol{x},t),\quad &&\Ddot{\Et}(\boldsymbol{x},t) &&= \Dot{\et}(\boldsymbol{x},t), \\
        \Dot{\pi}(\boldsymbol{x},t) &= p(\boldsymbol{x},t),\quad &&\Ddot{\pi}(\boldsymbol{x},t) &&= \Dot{p}(\boldsymbol{x},t).
    \end{alignat}
\end{subequations}
For the sake of simplicity, we introduce the following definitions:
\begin{subequations}
    \begin{alignat}{2}
        \operator f &= a\Ddot{f} + b\dot{f} + cf \\
        \tsigmaPML(\St,\pi) &= (\Dot{\St}-\alpha \Dot{\pi}\tidentity)\Lambdae + (\St-\alpha \pi\tidentity)\Lambdap
    \end{alignat}    
\end{subequations}
where $\operator$ denotes an operator that acts on any scalar, vector, or tensor function $f$ and $\tsigmaPML$ is the PML stress tensor. Thus, replacing the relations given in \eqref{eq: mixed variables} into \eqref{eq: time domain PML formulation} and introducing the previous definitions, the fully-mixed PML formulation becomes: find $\uu$, $\ww$, $\pi$, and $\St$ satisfying:
\begin{subequations}
    \label{eq: mixed time domain PML formulation}
    \begin{alignat}{2}
        \rho \operator \uu + \rho_f \operator \ww &= \nabla\cdot\tsigmaPML(\St,\pi) &&\quad\text{in }\Omega\times T \\ 
        \mathcal{D}(\operator \St) &= \frac{1}{2}\left\{ \nabla\uu\Lambdap + \Lambdap (\nabla\uu)^T + \nabla\Dot{\uu}\Lambdae + \Lambdae (\nabla\Dot{\uu})^T \right\} &&\quad\text{in }\Omega\times T \\
        -\nabla \Dot{\pi} &= \left(\Lambdaee\frac{\partial}{\partial t} + \Lambdapp\right)\left(\rho_f\Dot{\uu} + \rho_w\Dot{\ww} + \frac{\eta}{\kappa}\ww\right) &&\quad\text{in }\Omega\times T \\
        -\operator\pi &= \nabla\cdot\left\{M\left(\Lambdae\frac{\partial}{\partial t} + \Lambdap\right)(\alpha\uu + \ww)\right\} &&\quad\text{in }\Omega\times T
    \end{alignat}
\end{subequations}
where $\mathcal{D}$ denotes the compliance operator, which takes the stress tensor as argument and returns the strain tensor ($\et=\mathcal{D}(\tsigma)$). 

\textbf{Remark:} the tensor $\St$ is symmetric and, therefore, only three additional scalar functions are introduced as unknowns by the mixed PML formulation.


\section{Hybrid formulation}

With the fully-mixed formulation given in \eqref{eq: mixed time domain PML formulation}, the vector functions $\uu, \ww$, the scalars $p, \pi$, and the symmetric tensor field  $\St$ must be solved simultaneously on the whole domain $\Omega=\OmegaRD\cup\OmegaPML \subset \mathbb{R}^2$, which may lead to computationally expensive problems. Therefore, we define a hybrid formulation where the problem given in \eqref{eq: mixed time domain PML formulation} is split into two sub-problems defined separately on $\OmegaRD$ and $\OmegaPML$ but coupled through boundary conditions on $\GammaI$ (see Figure \ref{fig: PML reference domain}).

Let $\{\uu_1,\ww_1\}$ and $\{\uu_2,\ww_2\}$ be the solid and relative fluid displacements defined separately on $\OmegaRD$ and $\OmegaPML$, respectively. The hybrid PML formulation for the solid and fluid displacements, pore pressure, stress history, and pore pressure history reads: find $\{\uu_1,\ww_1,p\}$ and $\{\uu_2,\ww_2,\pi,\St\}$ satisfying:
\begin{subequations}
    \label{eq: hybrid Biot's-PML equations}
    \begin{alignat}{2}
        \rho\Ddot{\uu}_1 + \rho_f\Ddot{\ww}_1 &= \nabla\cdot\tsigma(\uu_1,p) &&\quad\text{in }\OmegaRD\times T \label{eq: hybrid Biot's-PML equations (a)} \\
        -\nabla p &= \rho_f\Ddot{\uu}_1 + \rho_w\Ddot{\ww}_1 + \frac{\eta}{\kappa}\Dot{\ww}_1 && \quad\text{in }\OmegaRD\times T \label{eq: hybrid Biot's-PML equations (b)} \\
        -\Dot{p} &= \nabla\cdot\left\{M(\alpha\Dot{\uu}_1 + \Dot{\ww}_1)\right\} && \quad\text{in }\OmegaRD\times T \label{eq: hybrid Biot's-PML equations (c)} \\
        \rho \operator \uu_2 + \rho_f \operator \ww_2 &= \nabla\cdot\tsigmaPML(\St,\pi) && \quad\text{in }\OmegaPML\times T \label{eq: hybrid Biot's-PML equations (d)} \\
        \mathcal{D}(\operator \St) &= \frac{1}{2}\left\{ \nabla\uu_2\Lambdap + \Lambdap (\nabla\uu_2)^T + \nabla\Dot{\uu}_2\Lambdae + \Lambdae (\nabla\Dot{\uu}_2)^T \right\} && \quad\text{in }\OmegaPML\times T \label{eq: hybrid Biot's-PML equations (e)} \\
        -\nabla \Dot{\pi} &= \left(\Lambdaee\frac{\partial}{\partial t} + \Lambdapp\right)\left(\rho_f\Dot{\uu}_2 + \rho_w\Dot{\ww}_2 + \frac{\eta}{\kappa}\ww_2\right) && \quad\text{in }\OmegaPML\times T \label{eq: hybrid Biot's-PML equations (f)} \\
        -\operator \pi &= \nabla\cdot\left\{M\left(\Lambdae\frac{\partial}{\partial t} + \Lambdap\right)(\alpha\uu_2 + \ww_2)\right\} && \quad\text{in }\OmegaPML\times T \label{eq: hybrid Biot's-PML equations (g)}        
    \end{alignat}
\end{subequations}
subject to zero initial values and the Dirichlet and Neumann boundary conditions listed below (see Figure \ref{fig: truncated semispace PML}).
\begin{subequations}
    \label{eq: hybrid boundary conditions}
    \begin{alignat}{2}
        \tsigma(\uu_1,p)\cdot\boldsymbol{n}_1 &= \boldsymbol{g} &&\quad\text{in }\Gamma_g\times T \label{eq: hybrid boundary conditions (a)} \\
        \tsigma(\uu_1,p)\cdot\boldsymbol{n}_1 = p\tidentity\cdot\boldsymbol{n}_1&=\boldsymbol{0} &&\quad\text{in }\Gamma^{\text{RD}}_N\times T \label{eq: hybrid boundary conditions (b)} \\
        \tsigmaPML(\St,\pi)\cdot\boldsymbol{n}_2 = \dot{\pi}\tidentity\cdot\boldsymbol{n}_2 &= \boldsymbol{0} &&\quad\text{in }\GammaNPML\times T \label{eq: hybrid boundary conditions (c)} \\
        \uu_2 &= \boldsymbol{0} &&\quad\text{in }\GammaDPML\times T \label{eq: hybrid boundary conditions (d)} \\
        \ww_2 &= \boldsymbol{0} &&\quad\text{in }\GammaDPML\times T \label{eq: hybrid boundary conditions (e)} \\
        \pi &= 0 &&\quad\text{in }\GammaDPML\times T \label{eq: hybrid boundary conditions (f)}
    \end{alignat}
\end{subequations}
where $\boldsymbol{n}_1$ and $\boldsymbol{n}_2$ are outward pointing normal vectors to $\OmegaRD$ and $\OmegaPML$ ($\boldsymbol{n}_1 = -\boldsymbol{n}_2$ in $\GammaI$), respectively (see Figure \ref{fig: PML truncation}).

Finally, to couple both equations, the continuity of displacements, tractions, and pressures must be imposed on the interface as follows:
\begin{subequations}
    \label{eq: interface conditions}
    \begin{alignat}{2}
        \tsigmaPML(\St,\pi)\boldsymbol{n}_1 + \tsigma(\uu_1,p)\boldsymbol{n}_2 &= 0 
        &&\quad\text{in }\GammaI\times T \label{eq: interface conditions (a)} \\
        \uu_1 &= \uu_2 &&\quad\text{in }\GammaI\times T \label{eq: interface conditions (b)} \\
        \ww_1 &= \ww_2 &&\quad\text{in }\GammaI\times T \label{eq: interface conditions (c)} \\
        p &= \dot{\pi} &&\quad\text{in }\GammaI\times T \label{eq: interface conditions (d)}
    \end{alignat}
\end{subequations}

\subsection{Extension to M-PML}
The previous formulations are prone to time instabilities depending of the media properties and the form of the stretching functions. To transform the PML layer in the fully-mixed and hybrid problems to the multiaxial case, the functions $\alpha$ and $\beta$ (see Equation \eqref{eq: PML function}) must be redefined as:
\begin{subequations}
    \begin{alignat}{2}
        \alpha(x) &= \alpha(y) = 1  \\
        \beta^*(x,y) &= \beta(x) + p^{(y/x)}\beta(y)\\
        \beta^*(x,y) &= \beta(y) + p^{(x/y)}\beta(x)
    \end{alignat}
\end{subequations}
where $p^{(y/x)}$ and $p^{(x/y)}$ are constant parameters that allow the fine-tuning of the M-PML layer. This modified definition of the scaling and attenuation functions does not introduce changes to the previous formulation.

\subsection{Variational formulation}

The interface conditions \eqref{eq: interface conditions (a)} and \eqref{eq: interface conditions (b)} are fulfilled by using a single continuous function space for $\uu_1$ and $\uu_2$. Similarly, \eqref{eq: interface conditions (c)} is fulfilled by using a single continuous function space for $\ww_1$ and $\ww_2$. A Lagrange multiplier is used for the condition \eqref{eq: interface conditions (d)}.

Thus, using the following function spaces
\begin{subequations}
    \label{eq: functional spaces}
    \begin{alignat}{2}
            V &= \{\uu\in[H^1(\Omega)]^2 ~\text{s.t.}~\uu=\boldsymbol{0} \text{ on } \GammaDPML\times T \} \\
            Q &= \{p\in L^2(\OmegaRD)\} \\
            Q_0 &= \{p\in L^2(\OmegaPML)~\text{s.t.}~p=0 \text{ on } \GammaDPML\times T\} \\
            T &= \{ S\in [L^2(\OmegaPML)]^{2\times 2}\} \\
            L &= \{ l\in L^2(\GammaI)\}
    \end{alignat}
\end{subequations}
the weak form of the system of PDEs given in \eqref{eq: hybrid Biot's-PML equations} reads: find $(\uu,\ww,p,\pi,\St,\lp)$ for all $(\tuu,\tww,\tp,\tpi,\tSt,\tlp)\in V\times V\times Q\times Q_0\times T\times L$ solution to:
\begin{subequations}
    \label{eq: variational formulation}
    \begin{alignat}{2}
        \innerRD{(\rho\Ddot{\uu} + \rho_w\Ddot{\ww})}{\tuu} + \tensorinnerRD{\tsigma(\uu,p)}{\nabla\tuu} &= \innerNRD{\boldsymbol{g}}{\tuu} \label{eq: variational formulation (a)}\\ 
        \innerPML{(\rho\operator\uu + \rho_f\operator\ww)}{\tuu}
        + \tensorinnerPML{\tsigmaPML(\St,\pi)}{\nabla\tuu} &= 0  \label{eq: variational formulation (b)}\\
        \innerRD{\left(\rho_f\Ddot{\uu} + \rho_w\Ddot{\ww} + \frac{\eta}{\kappa}\Dot{\ww}\right)}{\tww} - \scalarinnerRD{p}{\nabla\cdot\tww} + \innerg{p\boldsymbol{n}}{\tww} &= 0 \label{eq: variational formulation (c)}\\
        \innerPML{\left(\Lambdaee\frac{\partial}{\partial t} + \Lambdapp\right)\left(\rho_f\Dot{\uu} + \rho_w\Dot{\ww} + \frac{\eta}{\kappa}\ww\right)}{\tww} - \scalarinnerPML{\Dot{\pi}}{\nabla\cdot\tww} &= 0 \label{eq: variational formulation (d)}\\
        \scalarinnerRD{\Dot{p}}{\tp} + \scalarinnerRD{\nabla\cdot\left\{M(\alpha\Dot{\uu} + \Dot{\ww})\right\}}{\tp} &= 0 \label{eq: variational formulation (e)}\\
        \scalarinnerPML{\operator\pi}{\tpi} + \scalarinnerPML{\nabla\cdot\left\{M\left(\Lambdae\frac{\partial}{\partial t} + \Lambdap\right)(\alpha\uu + \ww)\right\}}{\tpi} &=0 \label{eq: variational formulation (f)}\\
        \tensorinnerPML{\mathcal{D}(\operator\St)}{\tSt} - \frac{1}{2}\tensorinnerPML{\{\nabla\uu\Lambdap + \Lambdap (\nabla\uu)^T + \nabla\Dot{\uu}\Lambdae + \Lambdae (\nabla\Dot{\uu})^T\}}{\tSt} &= 0 \label{eq: variational formulation (g)} \\
        \scalarinnerI{\tlp}{(p-\dot{\pi})} + \scalarinnerI{(\tp - \tpi)}{\lp} &= 0
    \end{alignat}
\end{subequations}
where $\lp$ is a Lagrange multiplier used to impose the coupling condition \eqref{eq: interface conditions (d)}.

\begin{figure}[h]
    \centering
    \subfloat[]{\includegraphics[width=0.45\textwidth]{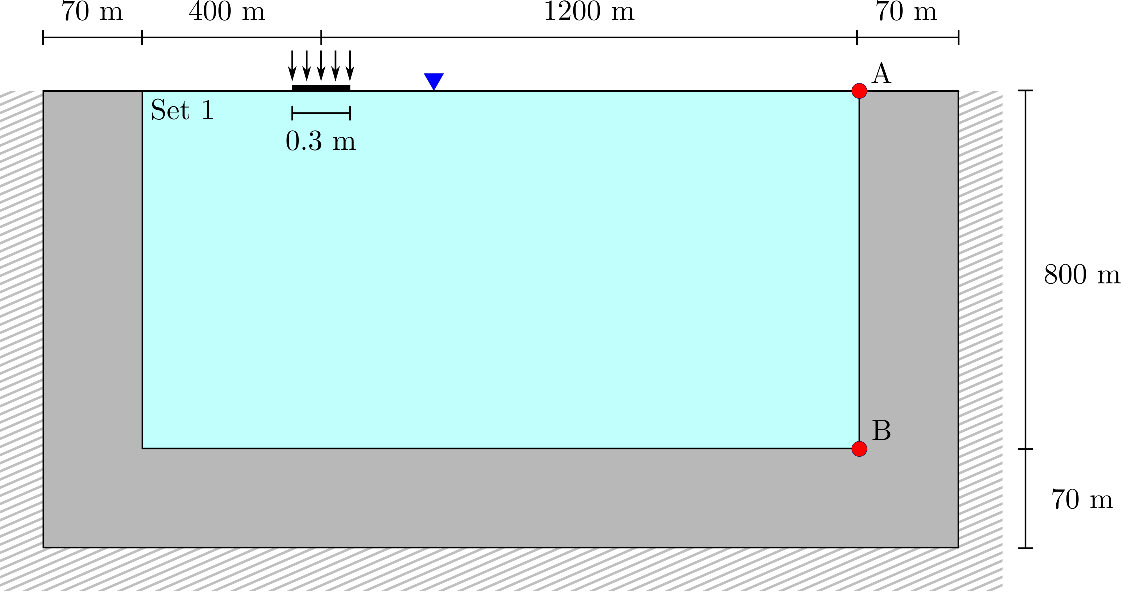} \label{fig: homogeneous experiment}}
    \hfil
    \subfloat[]{\includegraphics[width=0.45\textwidth]{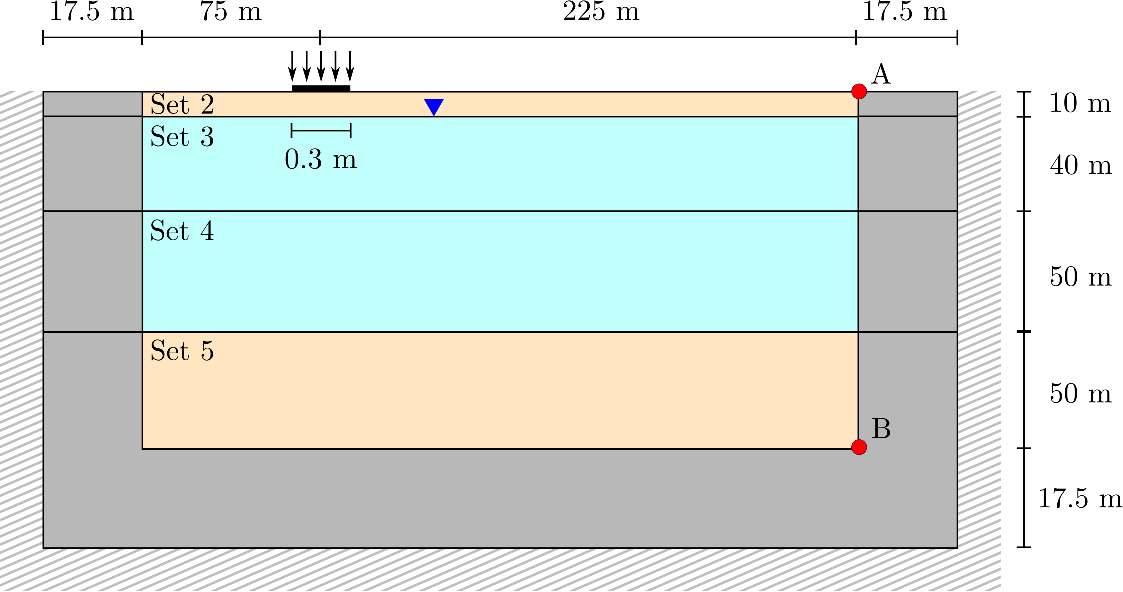} \label{fig: heterogeneous experiment}}
    
    \subfloat[]{\includegraphics[width=0.45\textwidth]{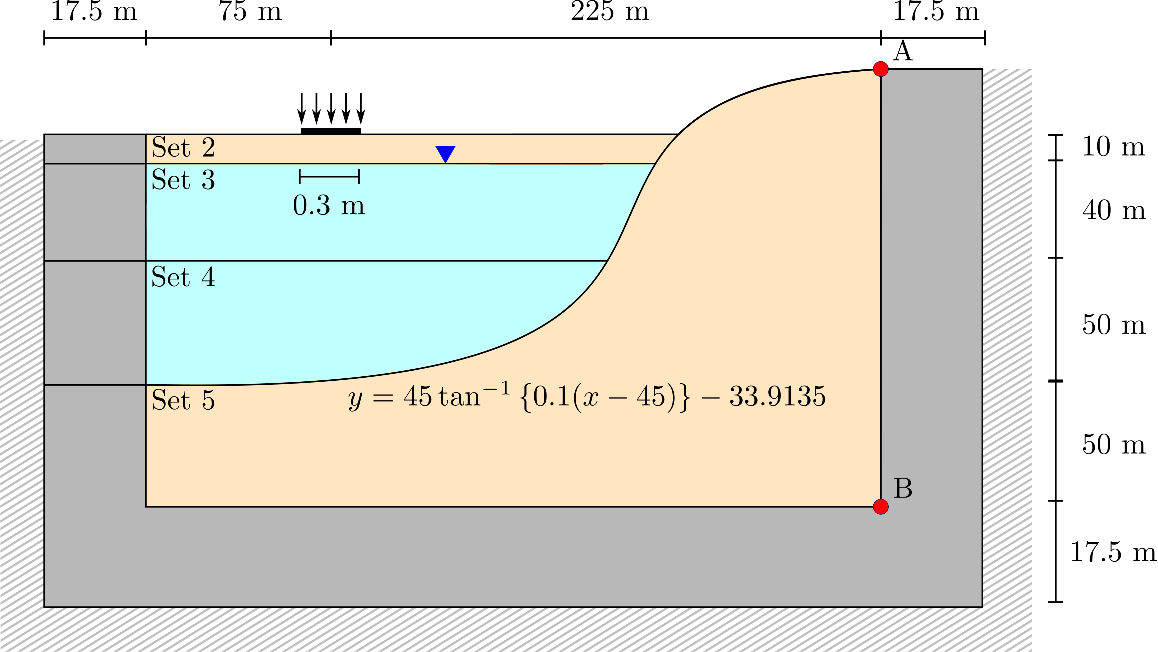}  \label{fig: complex heterogeneous experiment}}
    \subfloat[]{\includegraphics[width=0.45\textwidth]{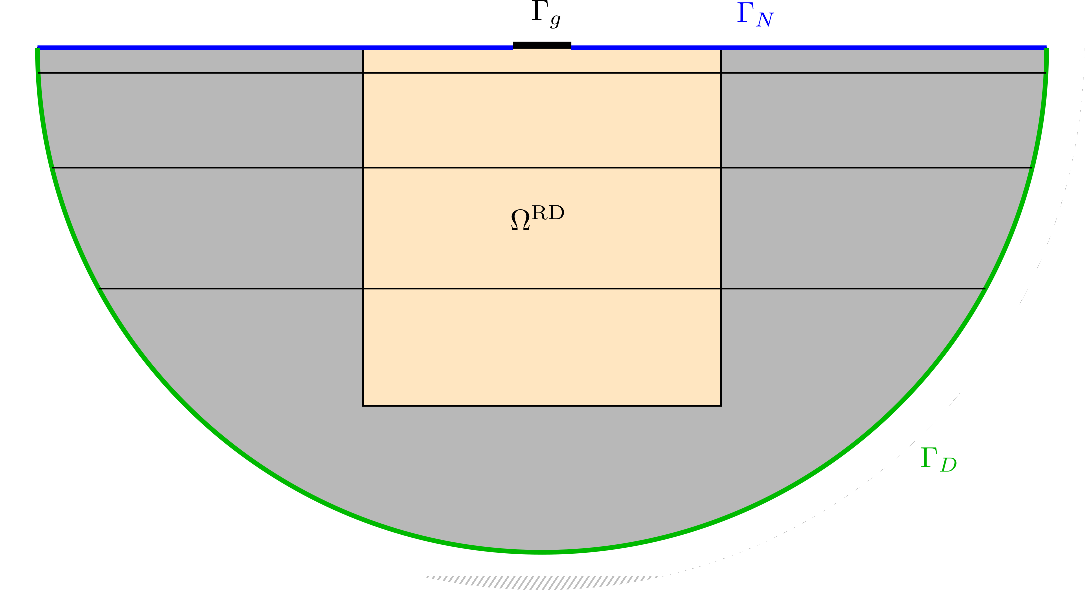} \label{fig: extended domain experiment}}
    
    \caption{Poroelastic domains used for the numerical experiments in (a) homogeneous, (b) horizontally layered, and (c) horizontally layered  with realistic interface to bedrock. The water table is marked with an inverted blue triangle. The red dots denote the locations where measurements of solid displacement ($\uu$), fluid velocity ($\ww$), and pressure ($p$) were taken. In Figure (d), a generic representation of the extended domains used in layered media is shown.}
    \label{fig: numerical experiments}
\end{figure}

\section{Numerical experiments}

Three experiments were developed to evaluate the performance and accuracy of the proposed fully-mixed and hybrid PML formulations. The first experiment (referred as Experiment 1) considers a homogeneous poroelastic half-space where the water table is located at the free surface. In the second (Experiment 2), a three horizontally layered media over a half-space and water table below the second layer was considered (Figure \ref{fig: heterogeneous experiment}). The third (Experiment 3) adds a more realistic stratification including an outcropping  (see Figure \ref{fig: complex heterogeneous experiment}). The material parameters used in the three cases are listed in Table \ref{tab1: physical parameters}. Set 1 corresponds approximately to a soft rock, while sets 2 to 5 represent standard soil parameters, from loose to dense sands. 

\begin{table}[]
\centering
\caption{Material parameters, characteristic frequency of the medium, and wave propagation velocities of the porous media considered in the experiments. The set 1 was taken from Tables 1 and 2 in \cite{DudleyWard2017}, whereas sets from 2 to 5 were defined by the authors to obtain soil-like wave velocities.}
\label{tab1: physical parameters}
\begin{tabular}{@{}lccccc@{}}
\toprule
                                & Set 1      & Set 2       & Set 3      & Set 4       & Set 5 \\ \midrule
$\rho_s~(\text{kg}/\text{m}^3)$ & $2650$     & $2600$      & $2600$     & $2600$      & $2600$ \\
$\rho_f~(\text{kg}/\text{m}^3)$ & $900$      & $1.29$      & $1000$     & $1000$      & $1.29$ \\
$K_s~(\text{N}/\text{m}^2)$     & $12\e{9}$  & $2.3\e{8}$  & $2.3\e{8}$ & $2.5\e{8}$  & $4.6\e{8}$  \\
$K_f~(\text{N}/\text{m}^2)$     & $2\e{9}$   & $1.4\e{5}$  & $2\e{9}$   & $2\e{9}$    & $1.4\e{5}$ \\
$K_b~(\text{N}/\text{m}^2)$     & $10\e{9}$  & $1.5\e{8}$  & $1.5\e{8}$ & $1.7\e{8}$  & $4.2\e{8}$ \\
$\mu_b~(\text{N}/\text{m}^2)$   & $5\e{9}$   & $1.33\e{8}$ & $1.33\e{8}$& $2.84\e{8}$ & $4.44\e{8}$ \\
$T$                             & $1.2$      & $1$         & $1$        & $1$         & $1$ \\
$\kappa~(\text{m}^2)$           & $1\e{-12}$ & $8\e{-9}$   & $8\e{-9}$  & $1\e{-10}$  & $1\e{-11}$ \\
$\phi$                          & $0.3$      & $0.2$       & $0.2$      & $0.2$       & $0.2$  \\
$\eta~(\text{Pa s})$            & $1\e{-3}$  & $2\e{-5}$   & $1\e{-3}$  & $1\e{-3}$   & $2\e{-5}$ \\ \midrule
$f_c~(\text{Hz})$ & 44210 & 62 & 4 & 318 & 49350 \\ 
$c_{s}~(\text{m/s})$ & 960 & 300 & 300 & 400 & 500 \\
$c_{1p}~(\text{m/s})$ & 2366 & 470 & 574 & 636 & 755 \\
$c_{2p}~(\text{m/s})$ & 775 & 329 & 425 & 513 & 330 \\ \bottomrule
\multicolumn{5}{l}{$f_c$: characteristic frequency of the medium, $c_{1p}$: fast primary wave velocity,} \\
\multicolumn{5}{l}{$c_{2p}$: slow primary wave velocity, $c_{s}$: shear wave velocity}
\end{tabular}
\end{table}

To obtain a reference solution, we solved Biot's equations given in \eqref{eq: Biot's equations} on an extended domain with dimensions large enough to avoid wave reflections at the exterior boundaries $\OmegaRD$. A zero displacement Dirichlet condition was imposed to both $\uu$ and $\ww$ in the exterior boundary $\Gamma_D$. A simulation time of 2.0 seconds (half of the runtime used for the PML experiments) was enough for comparison purposes in the extended domain simulations. For layered medium, the regular domain was embedded in the extended domain, and layers were extended to the exterior boundary (see Figure \ref{fig: extended domain experiment}). Additionally, a first-order paraxial boundary condition \cite{Akiyoshi1994} was implemented to solve Biot's equations also for comparison purposes. The PML domain shown in Figure \ref{fig: numerical experiments} was removed for paraxial simulations and replaced by this ABC at this boundary.

A vertical load defined by a Ricker wavelet was applied on a $0.3~\text{m}$ width stripe of the free surface in the three experiments (see Figure \ref{fig: homogeneous experiment}). The expression defining the source is as follows:
\begin{equation}
    \label{eq: Ricker wavelet}
    \boldsymbol{g}(t) = A\left[\begin{array}{c}
        0  \\
        S(t) 
    \end{array}\right] \quad \text{with } S(t)=\frac{(0.25u^2-0.5)e^{-0.25u^2}-13e^{13.5}}{0.5+13e^{13.5}} \quad \text{and }0\leq t \leq \frac{6\sqrt{6}}{\omega_r}
\end{equation}
where $A$ denotes the pulse amplitude and $u = \omega_r t - 3\sqrt{6}$. In the previous expression, $\omega_r=2\pi f_r$ is the characteristic central angular frequency of the pulse. In all the experiments $A=10^4~\text{N/m}$ and $f_r=15~\text{Hz}$ were used. The frequency of this source was adjusted to obtain a frequency spectrum similar to those obtained from geophysical seismic surveys, which makes it suitable for the simulations.

The porous media utilized in the simulations exhibit primarily dispersive behavior, indicated by the characteristic frequency of the media ($f_c$) satisfying the inequality $f_c=(\eta\phi)/(2\pi \rho_f T \kappa) > f_r=15~\text{Hz}$, where the slowest primary wave does not propagate \cite{biot_theory_1956,DudleyWard2017}. Although the medium generated by the physical parameters in Set 3 (see Table \ref{tab1: physical parameters}) does not display dispersive behavior, its shear wave velocity is slower than its slowest volumetric wave. Therefore, the discretization parameters were chosen by considering only the fastest primary wave and the shear wave velocities for all media. The domain was discretized using triangular cells, and the element size ($\Delta x$) was adjusted to have a minimum of 12 elements per shortest wavelength, with the biggest possible element being chosen. In all simulations, elements in the vicinity of the surface load were refined to a size of $\Delta x=0.15~\text{m}$.

A third-order version of the scaling and attenuation profiles given in \eqref{eq: PML polynomial profiles} (with $m=3$) was used for the simulations. The width of the PML, denoted by $L_{\text{PML}}$, was chosen to be ten times the element size once $\Delta x$ was fixed (see Table \ref{tab1: physical parameters}). The constant $\beta_0$ was calculated using the expression \eqref{eq: PML scaling constants} with $R=10^{-4}$, while $\alpha_0$ was fixed at 5. For all the experiments considering M-PML stretching functions, $p^{(y/x)}$ and $p^{(x/y)}$ were used as $0.01$. A summary of the discretization and PML parameters is presented in Table \ref{tab2: summary of parameters}.

The Newmark-$\beta$ method with $\beta=1/4$ and $\gamma=1/2$ (i.e., without numerical damping) was used for time discretization \cite{newmark_method_1959}. The time-step $\Delta t$ was calculated using the Courant-Friedrichs-Lewy criteria:
\begin{equation}
    \Delta t < \text{CFL}\frac{\Delta x}{c_{1p}}
\end{equation}
where $c_{1p}$ represents the velocity of the fast primary wave and $\text{CFL}=0.75$ is the Courant-Friedrichs-Lewy number.

\begin{table}[t]
\caption{Element sizes and PML parameters used in the three experiments with homogeneous and layered media}
\label{tab2: summary of parameters}
\centering
\begin{tabular}{@{}lccc|ccc@{}}
\toprule
              & \multicolumn{3}{c}{General parameters}                                  & \multicolumn{3}{c}{PML parameters}                       \\ \midrule
Experiment        & $\Delta x_{\text{global}}~\text{(m)}$ & $\Delta x_{\text{source}}~\text{(m)}$ & $\Delta t~\text{(s)}$ & $R$       & $L_{\text{PML}}~\text{(m)}$ & $\alpha_0$ \\ \midrule
1   & $7.8$                   & $0.15$                 & $10^{-3}$              & $10^{-4}$ & $78$                        & $5$            \\
2 and 3 & $1.4$                & $0.15$                 & $10^{-3}$              & $10^{-4}$ & $14$                      & $5$            \\ \bottomrule
\multicolumn{6}{l}{$\Delta x_{\text{global}}$: global element size; $\Delta x_{\text{local}}$: element size refinement near to external source} \\
\multicolumn{6}{l}{$\Delta t$: time step; $R$: reflection coefficient; $L_{\text{PML}}$: width of the PML layer}
\end{tabular}
\end{table}

\subsection{Metrics for performance evaluation}

To evaluate the performance of the fully-mixed PML, hybrid PML, and paraxial methods, the poroelastic energy on $\OmegaRD$ was estimated and compared against the reference solution using the following expression:
\begin{alignat}{2}
    \label{eq: energy}
    E(t_k) &= \frac{1}{2}\innerRD{\rho\Dot{\uu}(\boldsymbol{x},t_k)}{\Dot{\uu}(\boldsymbol{x},t_k)} + \frac{1}{2}\int_{\OmegaRD}\rho C\et(\uu(\boldsymbol{x},t_k),t_k) : \et(\uu(\boldsymbol{x},t_k),t_k)~d\Omega \nonumber \\
    &+ \frac{1}{2}\innerRD{\rho_w \Dot{\ww}(\boldsymbol{x},t_k)}{\Dot{\ww}(\boldsymbol{x},t_k)} + \frac{1}{2}\innerRD{\frac{1}{M}p(\boldsymbol{x},t_k)}{p(\boldsymbol{x},t_k)} \\
    &+ \innerRD{\rho_f \Dot{\uu}(\boldsymbol{x},t_k)}{\Dot{\ww}(\boldsymbol{x},t_k)} \nonumber
\end{alignat}

Finally, we obtained traces of $\uu$, $\ww$, and $p$ at different locations $\boldsymbol{x}_i$ in $\OmegaRD$ (see Figure \ref{fig: numerical experiments}). Normalized error metric are defined as:
\begin{subequations}
    \label{eq: trace error}
    \begin{alignat}{2}
        e_{\uu}(\boldsymbol{x}_i,t_k) &= \frac{\lVert \uu_{\text{ref}}(\boldsymbol{x}_i,t_k) - \uu(\boldsymbol{x}_i,t_k) \rVert_2}{\max_{t_k}\lVert\uu_{\text{ref}}(\boldsymbol{x}_i,t_k)\rVert_2} \\
        e_{\ww}(\boldsymbol{x}_i,t_k) &= \frac{\lVert \ww_{\text{ref}}(\boldsymbol{x}_i,t_k) - \ww(\boldsymbol{x}_i,t_k) \rVert_2}{\max_{t_k}\lVert\ww_{\text{ref}}(\boldsymbol{x}_i,t_k)\rVert_2} \\
        e_{p}(\boldsymbol{x}_i,t_k) &= \frac{\lvert p_{\text{ref}}(\boldsymbol{x}_i,t_k) - p(\boldsymbol{x}_i,t_k) \rvert}{\max_{t_k}\lvert p_{\text{ref}}(\boldsymbol{x}_i,t_k) \rvert}
    \end{alignat}
\end{subequations}
where $\lvert \cdot \rvert$ is the absolute value and $\lVert \boldsymbol{f} \rVert_2 = \sqrt{\sum_i f_i^2}$ the vector 2-norm. The sub-index $()_{\text{ref}}$ denotes the reference solution obtained in the extended domain simulations.

\subsection{Implementation}

All experiments were solved using the open-source computing platform FEniCS \cite{Alnaes2015,logg_automated_2012}. To implement the hybrid PML problem, Multiphenics \cite{multiphenics} was used as a complementary tool. The finite element meshes were generated using the Frontal-Delaunay algorithm in Gmsh \cite{geuzaine_gmsh_2009}. Discretization of the displacements ($\uu$ and $\ww$) and pressures ($p$ and $\pi$) was carried out using continuous Lagrange polynomials of second and first order, respectively. The stress history $\St$ was discretized using discontinuous Lagrange polynomials of first order.

The Multifrontal Massively Parallel sparse direct Solver (MUMPS) and the Generalized Minimal Residual Method (GMRES) solvers were used to solve the linear systems obtained after assembling the discrete weak forms of the problems \cite{petsc-web-page,petsc-efficient,petsc-user-ref}. FEniCS is built with PETSc as linear algebra backend \cite{petsc-web-page,petsc-efficient,petsc-user-ref} and supports both solvers by default. The iterative solver was used only for the extended domain simulations with relative and absolute tolerances of $10^{-7}$ and $10^{-9}$, respectively \cite{petsc-web-page}. To accelerate the convergence, the linear system of the extended domain simulation was right-preconditioned using the Parallel ILU preconditioner HYPRE-Euclid \cite{hysom_scalable_2001}. The direct solver was used with default parameters \cite{Alnaes2015}.

\section{Results}
In the upcoming sections, we will present and analyze energy graphs, traces, error traces, and snapshots of the propagating waves for all three experiments. Through these analyses, we aim to provide a comprehensive understanding of the simulations and their outcomes.

\subsection{Experiment 1: homogeneous half-space}

\begin{figure}[h!]
    \centering
    \subfloat[]{\includegraphics[width=0.49\textwidth]{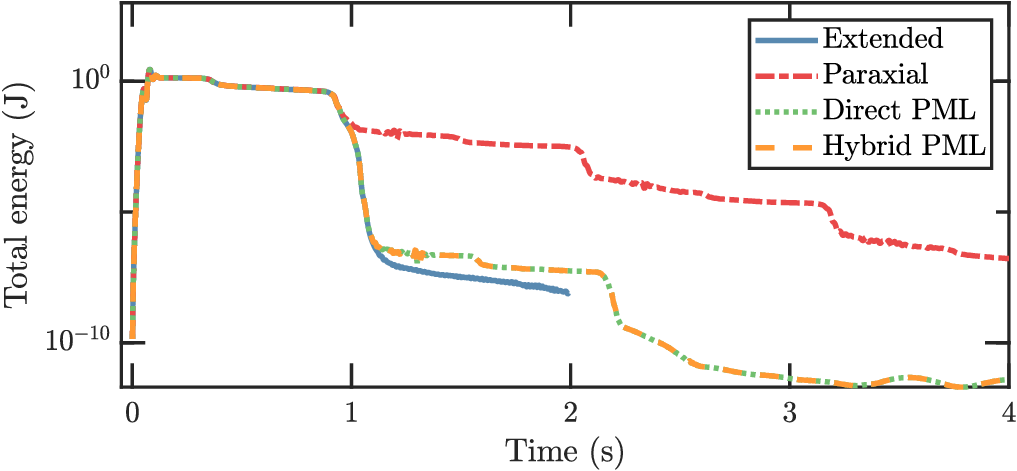}\label{fig: exp 1 energy a}} \hfill
    \subfloat[]{\includegraphics[width=0.49\textwidth]{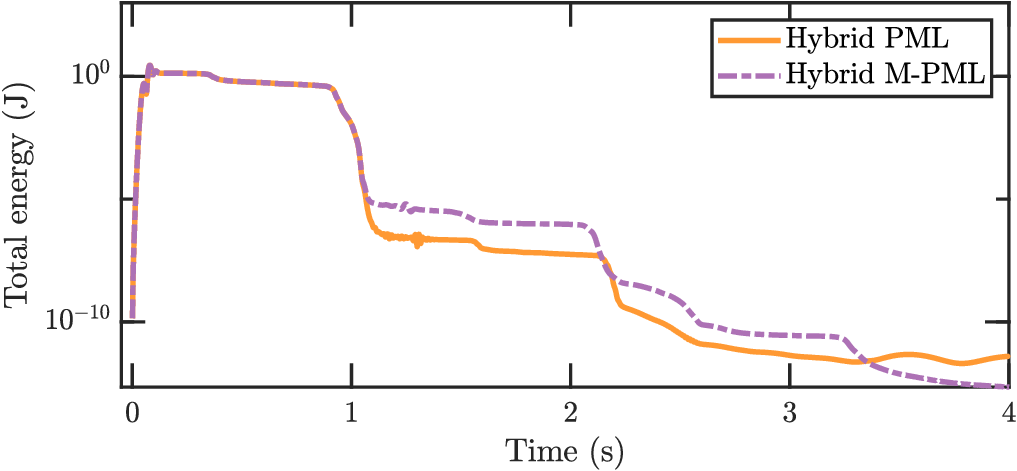}\label{fig: exp 1 energy b}}
    \caption{Energies estimated on $\OmegaRD$ using Equation \eqref{eq: energy} for the horizontally layered experiment. (a) Results show the extended, paraxial, fully-mixed, and hybrid PML results, and (b) a comparison between the hybrid PML and hybrid M-PML simulations. Both fully-mixed and hybrid formulations give identical results.}
    \label{fig: exp 1 energy}
\end{figure}

Figure \ref{fig: exp 1 energy a} shows the poroelastic energy estimated using \eqref{eq: energy} for extended, paraxial, and PML simulations in a homogeneous half-space domain. The results demonstrate good agreement between the PML and reference solutions in the first 2 seconds of simulation time, although small differences are observed. The energy decay rate of the PML simulations is greater than that of the paraxial case, which consistently decays but a lower rate. The energy obtained from the hybrid and fully-mixed PML formulations showed no observable differences, indicating that both methods provide equivalent solutions (see Figure \ref{fig: exp 1 energy a}). Additionally, the number of degrees of freedom (DOFs) in the hybrid case is approximately 1.7 times less than the fully-mixed problem (as shown in Table \ref{tab3: degrees of freedom}) because the tensor $\St$ does not need to be solved in $\OmegaRD$. Consequently, the hybrid formulation is significantly less computationally expensive than the fully-mixed form while maintaining the same properties.

\begin{table}[h!]
    \centering
    \caption{Number of degrees of freedom (DOFs) solved for extended, paraxial, fully-mixed PML, and hybrid PML formulations.}
    \label{tab3: degrees of freedom}
    \begin{tabular}{@{}lrrr@{}}
        \toprule
                        & Experiment 1 & Experiment 2 & Experiment 3 \\ \midrule
        Extended        & 8,550,901      & 13,043,765     & 12,759,950     \\
        Paraxial        & 425,518       & 462,513       & 480,229       \\
        Fully-mixed PML & 1,054,006      & 1,137,430      & 1,181,065      \\
        Hybrid PML      & 607,424       & 651,758       & 676,810       \\ \bottomrule
    \end{tabular}
\end{table}

Figure \ref{fig: exp 1 energy b} compares the energies obtained from the hybrid PML and M-PML formulations. During the first second of the simulation, the results are similar, and only small differences are observed. As the simulation progresses, the energy obtained with M-PML stretching functions is slightly larger than that obtained with uniaxial functions, indicating slightly worse performance. However, during the last second of simulation, the energy of the uniaxial case stops decaying showing even a slightly increase, while the energy of the M-PML case continues decaying. The reduced performance of M-PML (compared to PML) during the first seconds of simulation is because the stretching functions are not perfectly matched in $\GammaI$, generating spurious reflections. In fact, a M-PML absorbing layer can be interpreted as a sponge rather than a PML \cite{martin_high-order_2010}. This is because the coupling of two damping directions causes the loss of the perfectly matched layer characteristic of Berenger's technique \cite{Berenger1994}. Thus, the theoretical reflection coefficient for an infinite M-PML is not longer zero prior the discretization.

\begin{figure}[h!]
    \centering
    \subfloat[A Point]{\includegraphics[width=0.49\textwidth]{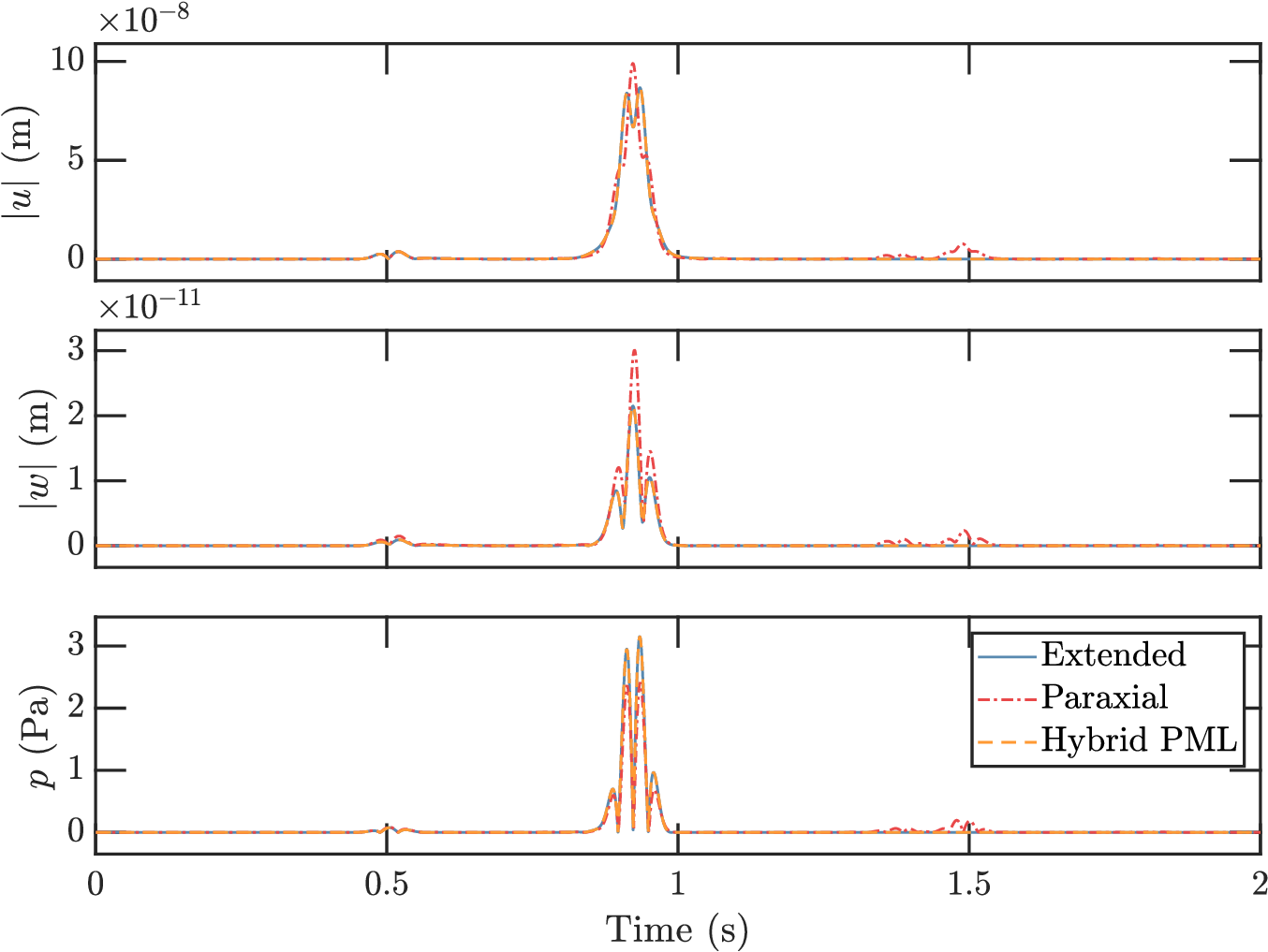}}\hfil
    \subfloat[B Point]{\includegraphics[width=0.49\textwidth]{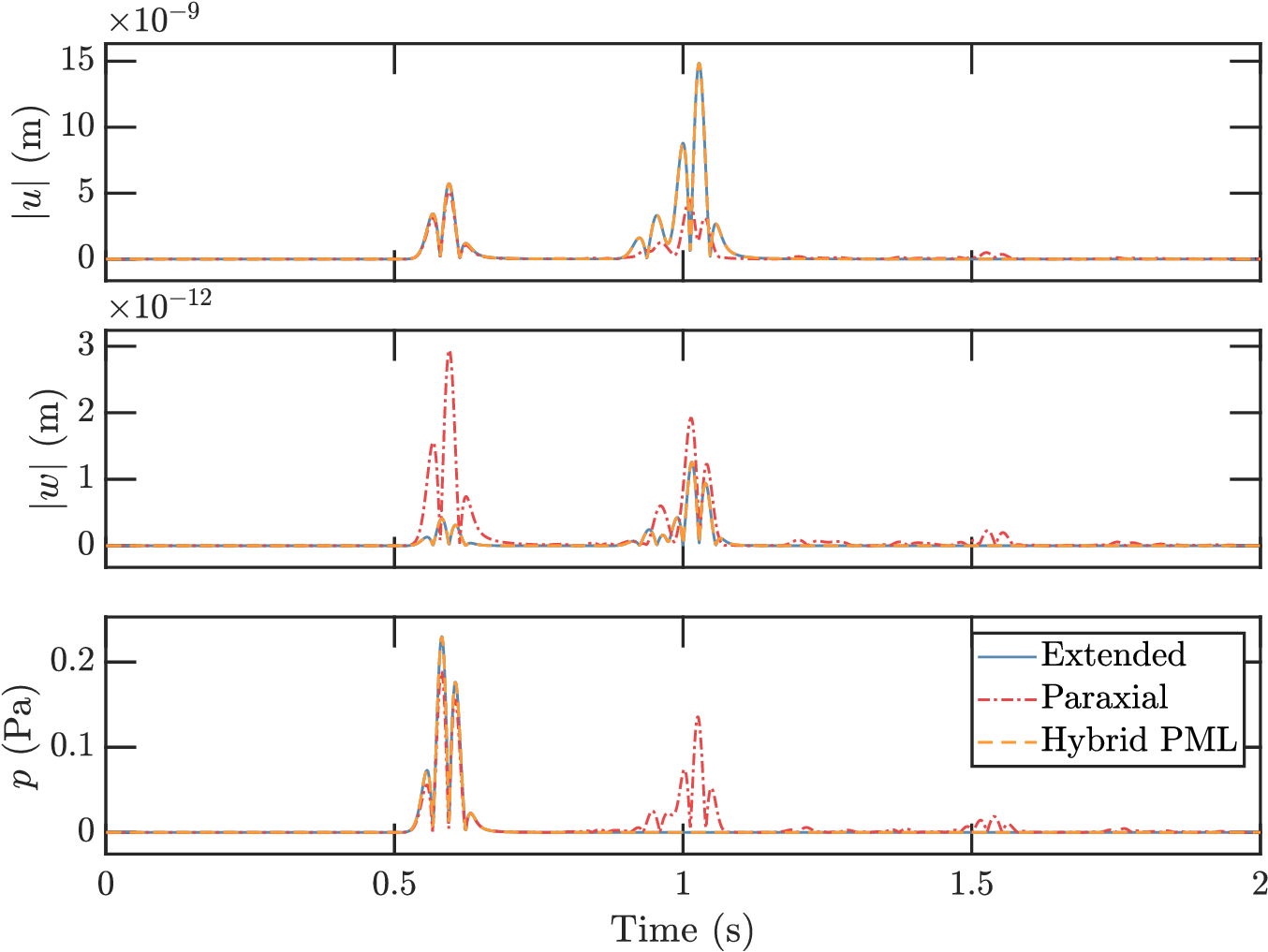}}

    \caption{Traces of $\uu$, $\ww$, and $p$ for Experiment 1 at the points highlighted in Figure \ref{fig: homogeneous experiment}. Results obtained using the hybrid PML formulation show good agreement with the reference solution and no spurious reflections are observed.}
    \label{fig: traces 1st experiment}
\end{figure}

\begin{figure}[h!]
    \centering
    \subfloat[Location A]{\includegraphics[width=0.49\textwidth]{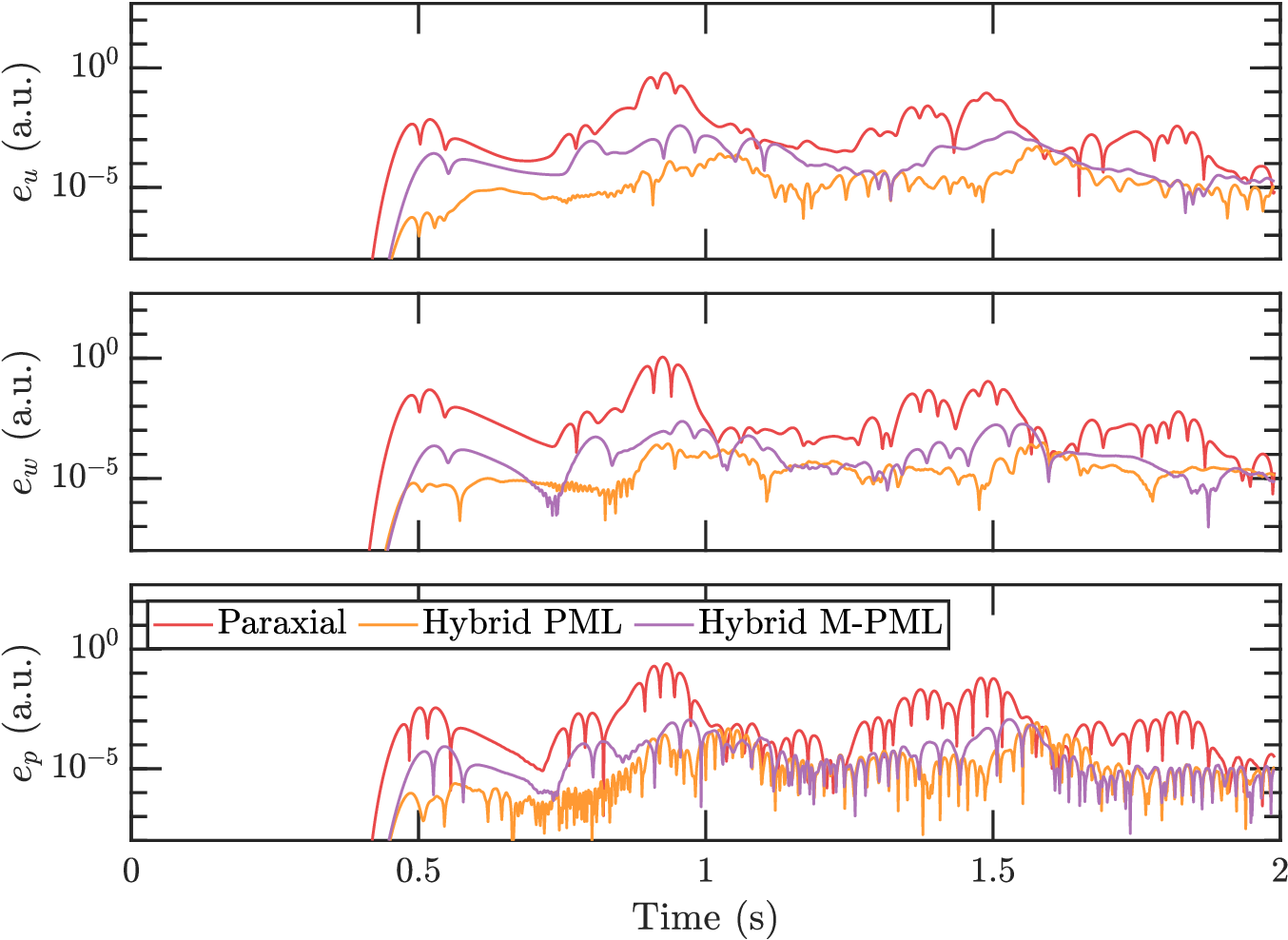}} \hfil
    \subfloat[Location B]{\includegraphics[width=0.49\textwidth]{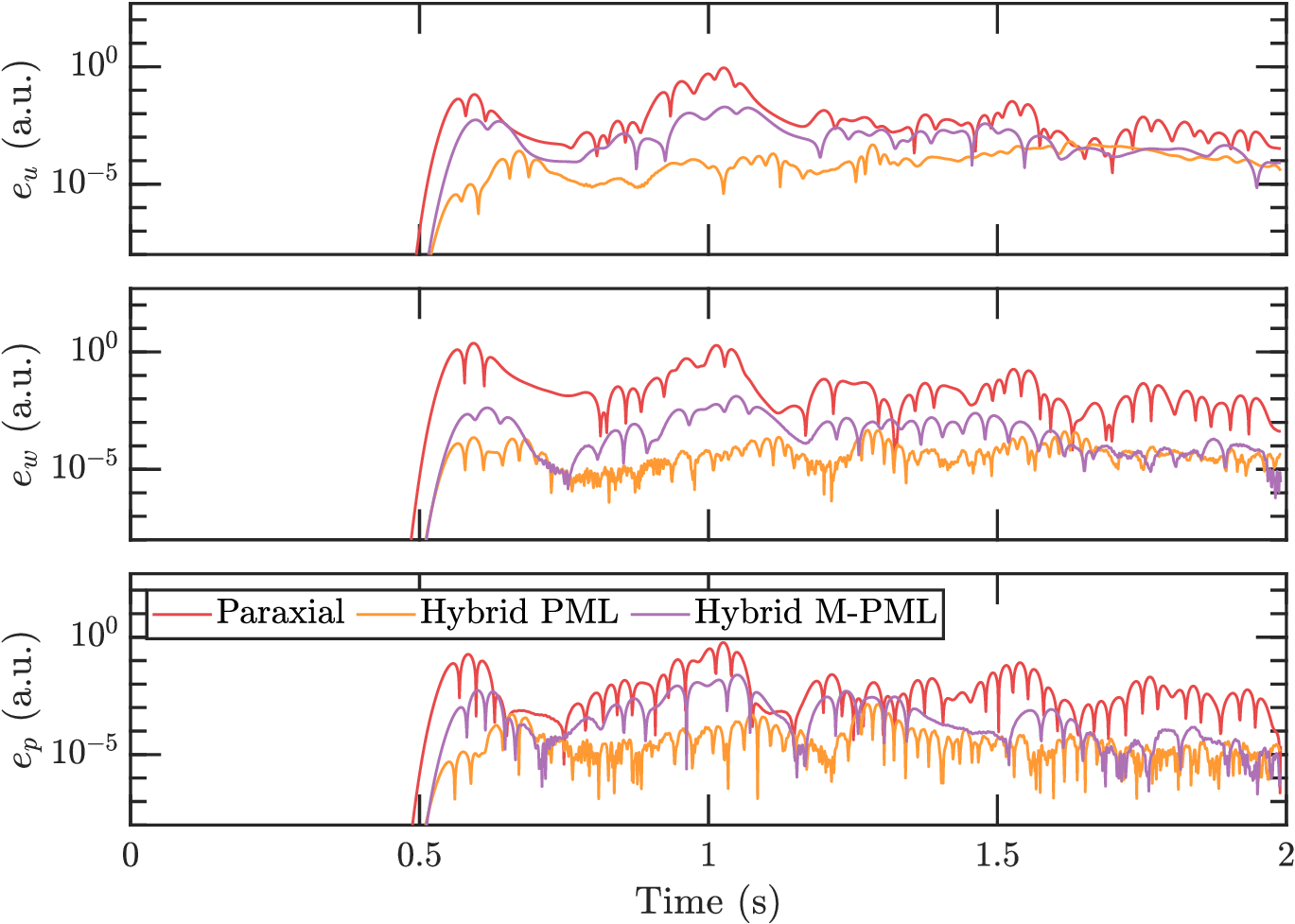}}

    \caption{Errors in the traces of $\uu$, $\ww$, and $p$ estimated using \eqref{eq: trace error} for the Experiment 1 at the highlighted locations in Figure \ref{fig: homogeneous experiment}. The vertical axis is logarithmic to facilitate visualization of differences. At the beginning of the simulation, all formulations yielded errors close to zero, which are omitted in the plot. As the simulation progressed, errors obtained using the hybrid PML and M-PML are smaller than those obtained using the paraxial case, although M-PML showed slightly worse performance between 0.5 and 1.5 s approximately because the stretching functions are not perfectly matched in this case. \ref{fig: exp 1 energy b}.}
    \label{fig: error of traces 1st experiment}
\end{figure}

\begin{figure}[h!]
    \centering
    \subfloat[$t=0.2$ s]{\includegraphics[scale=0.15]{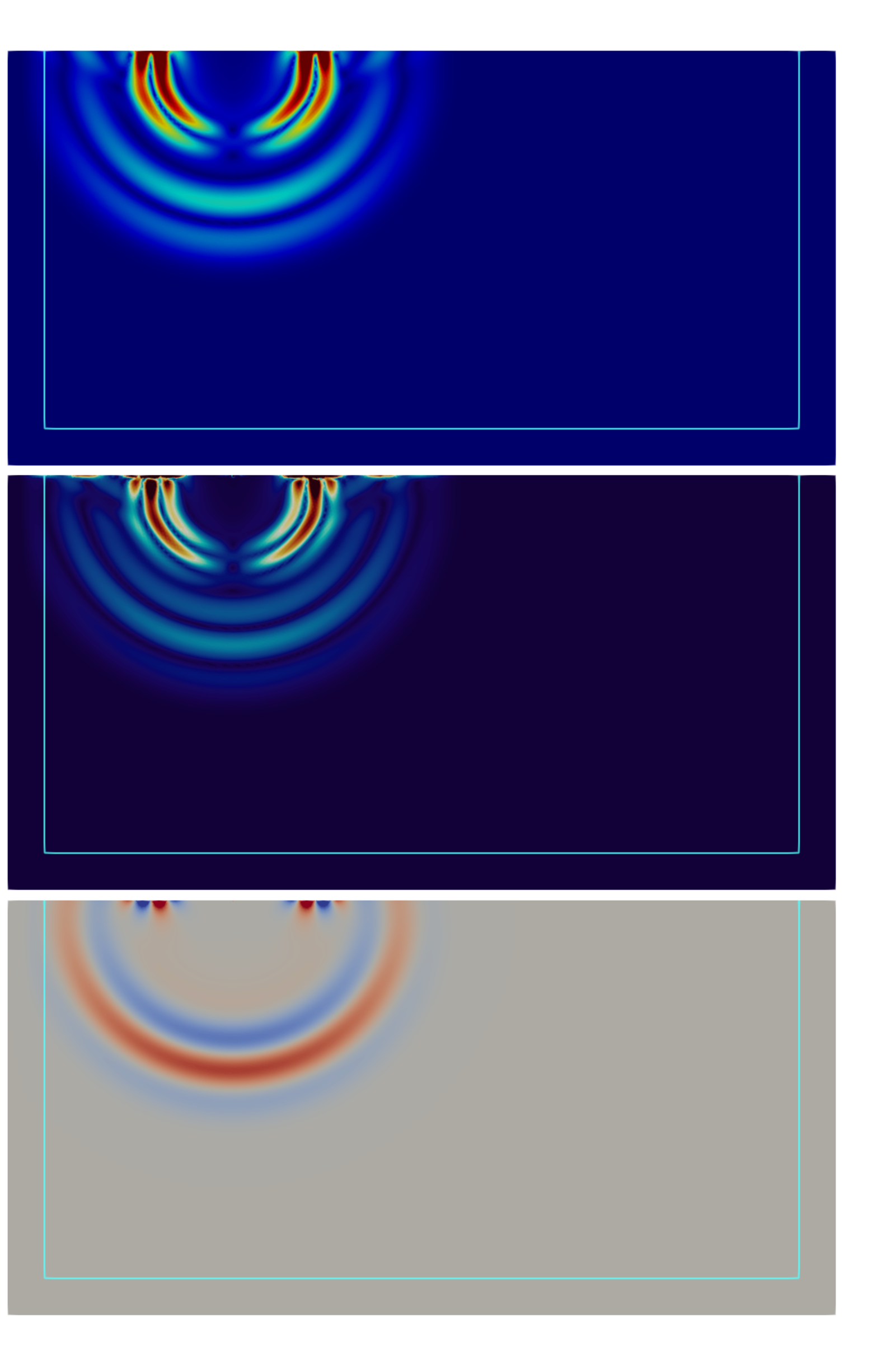} \label{fig: Exp 1 screenshots (a)}}
    \subfloat[$t=0.5$ s]{\includegraphics[scale=0.15]{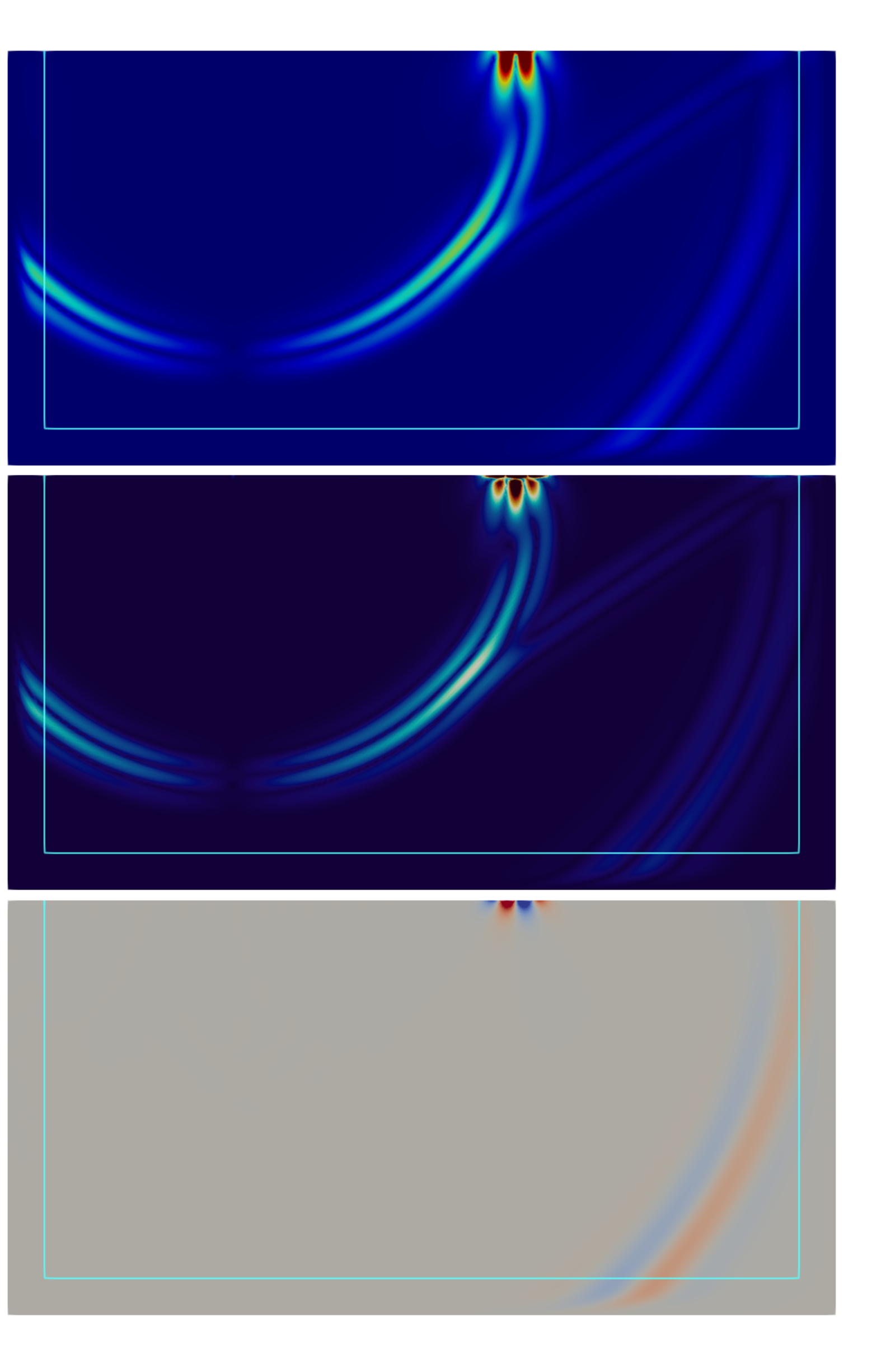} \label{fig: Exp 1 screenshots (b)}}
    \subfloat[$t=1$ s]{\includegraphics[scale=0.15]{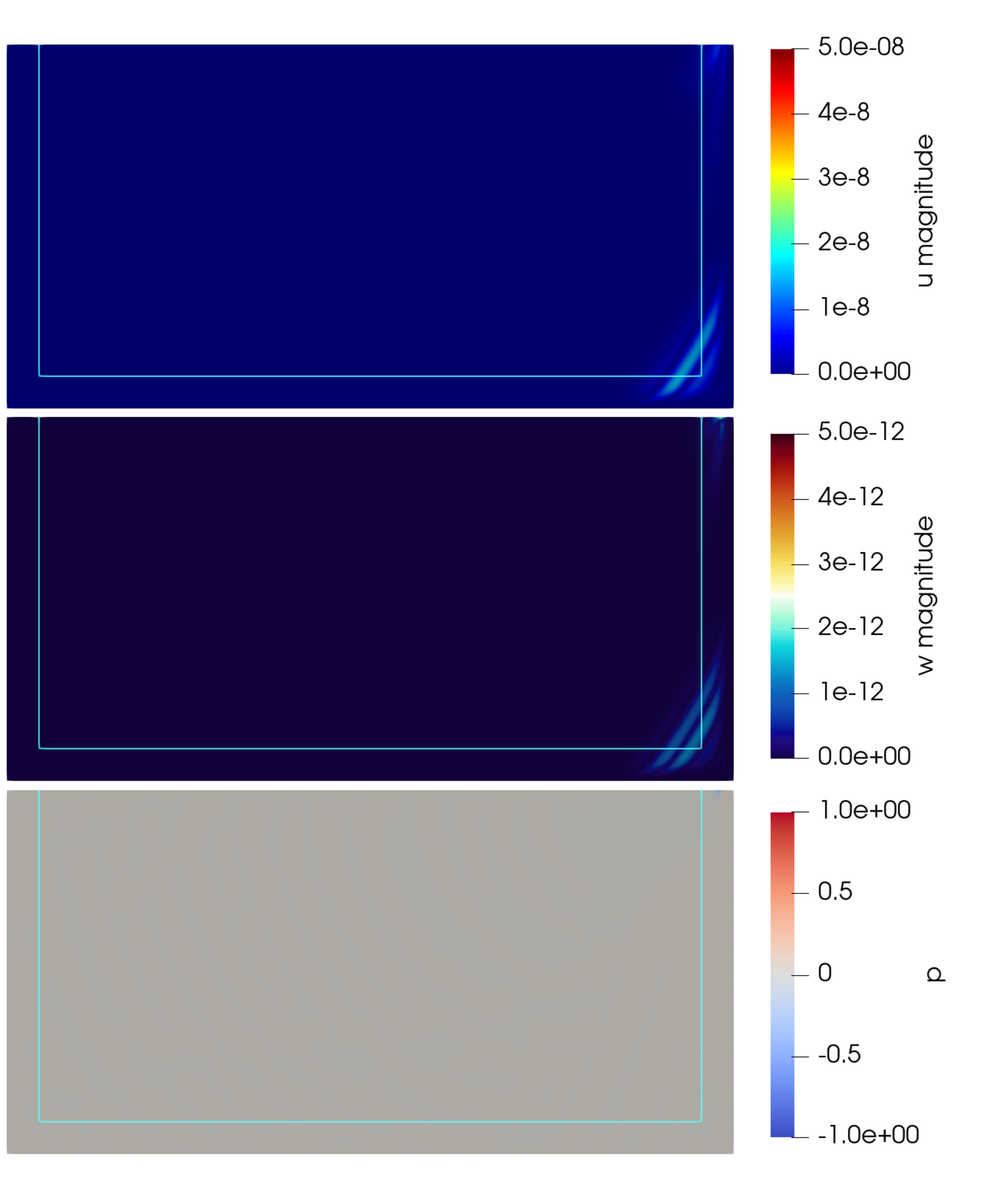} \label{fig: Exp 1 screenshots (c)}}

    \caption{Screenshots of the hybrid PML simulations at different time steps. The first row in the figure shows the solid displacement ($\uu$), the second the relative fluid displacement ($\ww$), and the third row the fluid pressure ($p$). No evident spurious reflections are observed at the interface between $\OmegaRD$ and $\OmegaPML$.}
    \label{fig: Exp 1 screenshots}
\end{figure}

Traces of the solutions and errors calculated using \eqref{eq: trace error} for the two locations shown in Figure \ref{fig: homogeneous experiment} are presented in Figures \ref{fig: traces 1st experiment} and \ref{fig: error of traces 1st experiment}. There is a good match between the hybrid PML and the extended domain simulations at both locations, in contrast to the paraxial case where reflections are observed. Looking at the errors, the superior performance of the hybrid PML method is evident, showing an improvement of at least three orders of magnitude compared to the paraxial case. In the same figure, the error obtained by using M-PML stretching functions in the hybrid simulation is depicted. As mentioned previously, results obtained using M-PML compared to PML are slightly worse between 0.5 and 1.5 s approximately because the stretching functions are not perfectly matched in this case. However, the results are still better than those obtained with the paraxial boundary conditions and more stable in time compared to the hybrid PML simulation (see Figure \ref{fig: exp 1 energy}).  Finally, screenshots of the solutions at different times can be found in Figure \ref{fig: Exp 1 screenshots}.

\subsection{Experiment 2: horizontally-layered domain}

\begin{figure}[h!]
    \centering
    \subfloat[]{\includegraphics[width=0.49\textwidth]{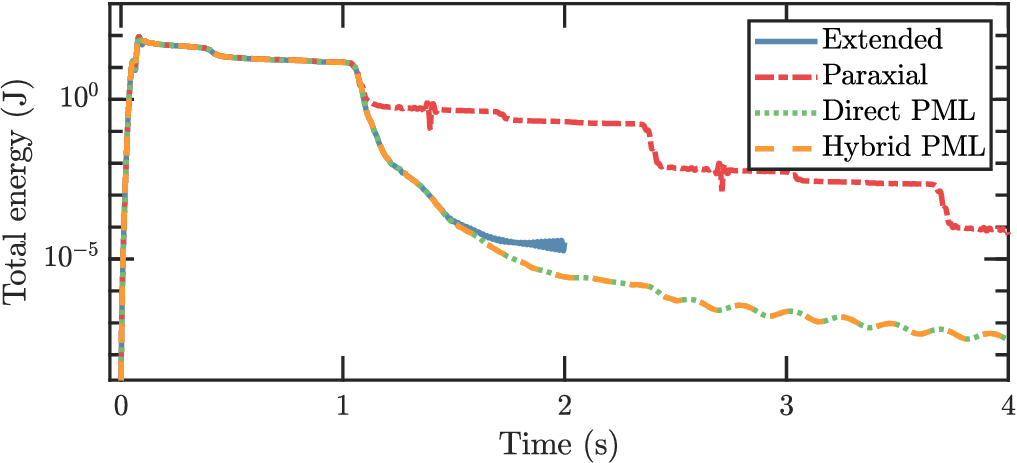} \label{fig: heterogeneous energy a}} \hfill
    \subfloat[]{\includegraphics[width=0.49\textwidth]{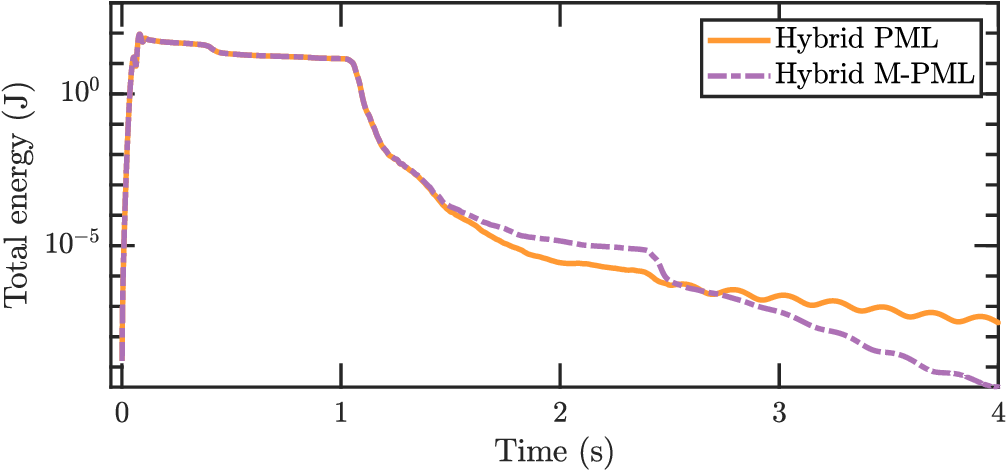} \label{fig: heterogeneous energy b}}
    \caption{Energies estimated on $\OmegaRD$ using Equation \eqref{eq: energy} for the horizontally layered experiment. (a) Results show the extended, paraxial, fully-mixed, and hybrid PML results, and (b) a comparison between the hybrid PML and hybrid M-PML simulations. Both fully-mixed and hybrid formulations give identical results.}
    \label{fig: heterogeneous energy}
\end{figure}

\begin{figure}[h!]
    \centering
    \subfloat[$t=0.3$ s]{\includegraphics[scale=0.15]{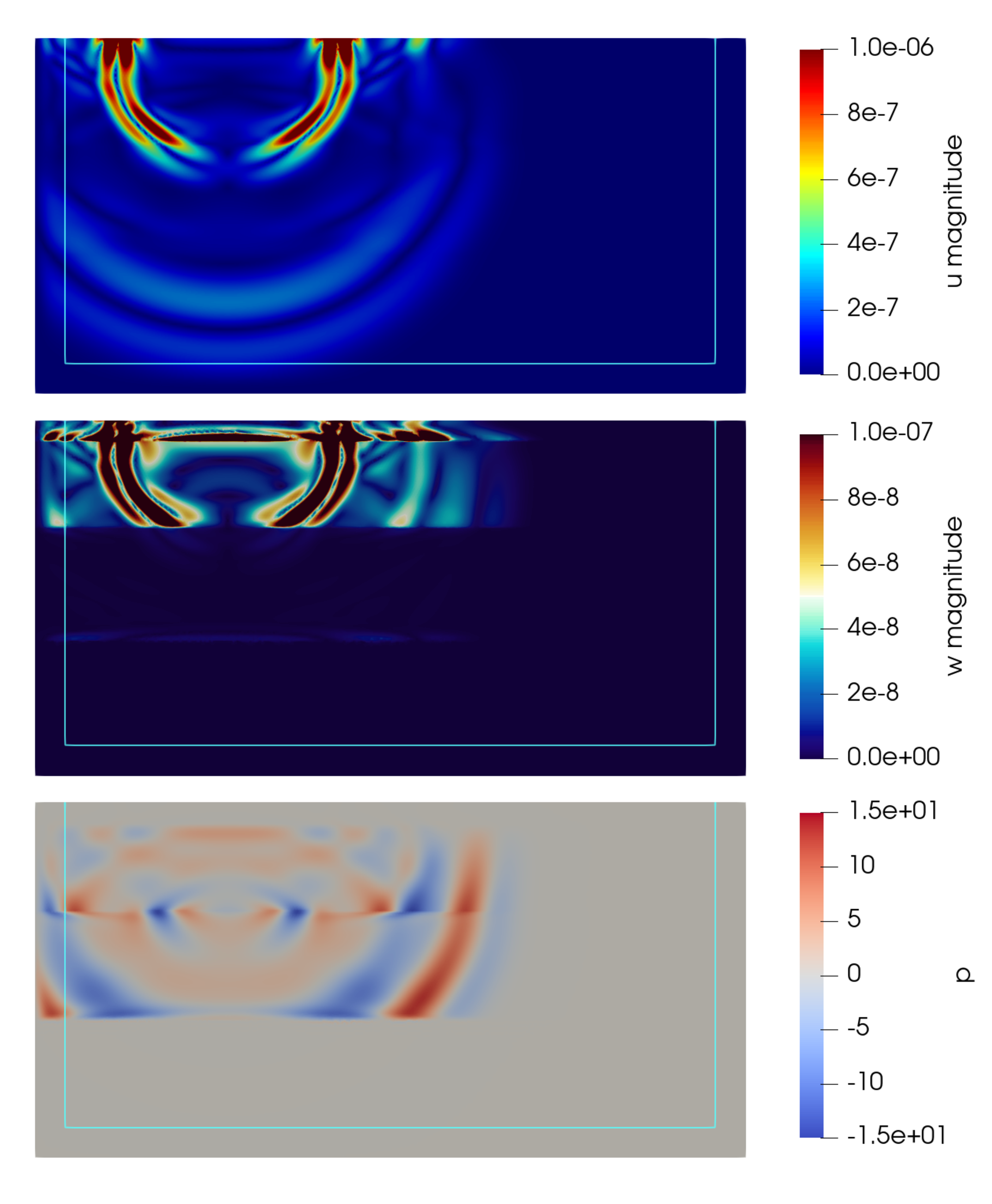} \label{fig: Exp 2 screenshots (a)}}
    \subfloat[$t=0.5$ s]{\includegraphics[scale=0.15]{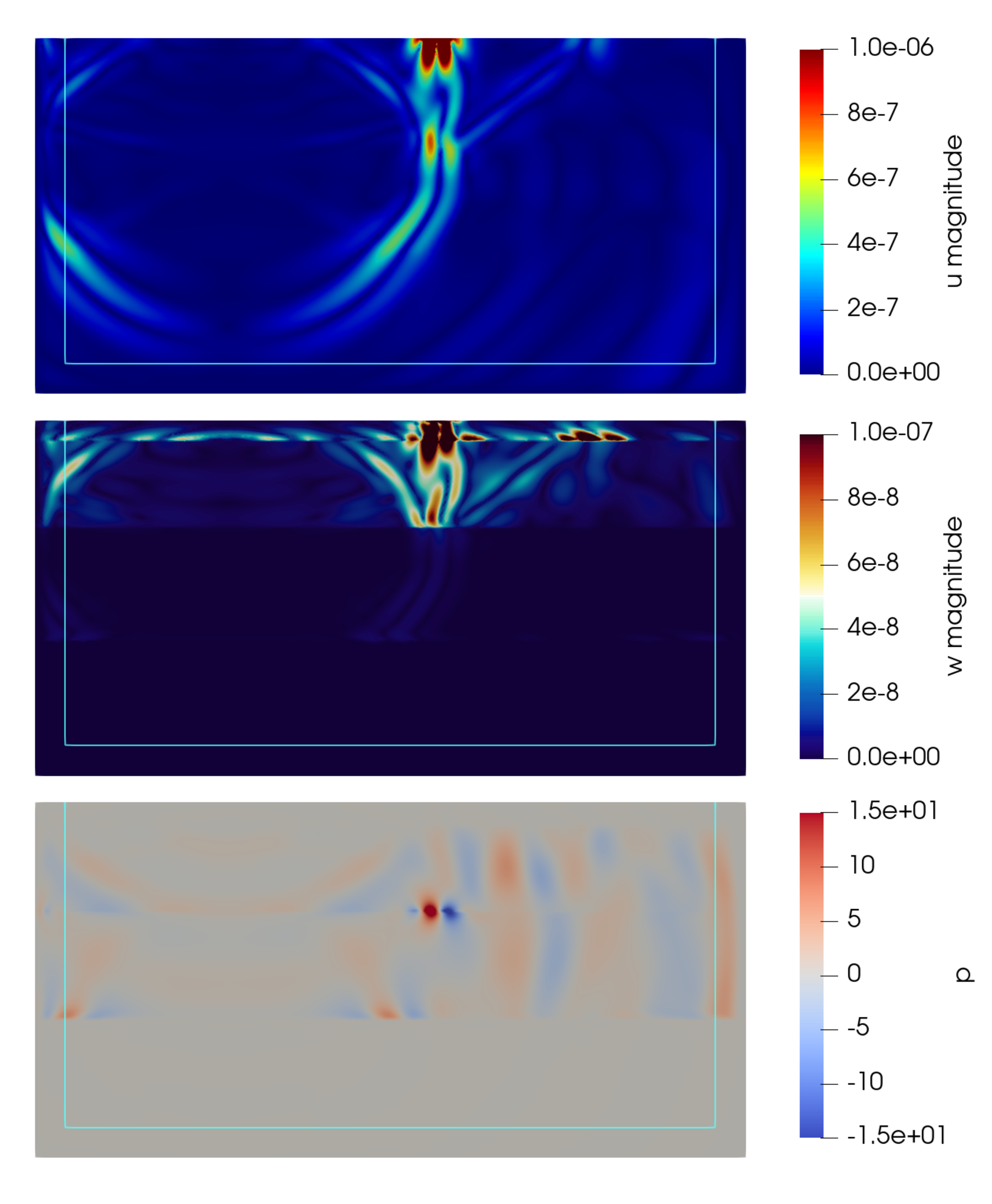} \label{fig: Exp 2 screenshots (b)}}
    \subfloat[$t=1.2$ s]{\includegraphics[scale=0.15]{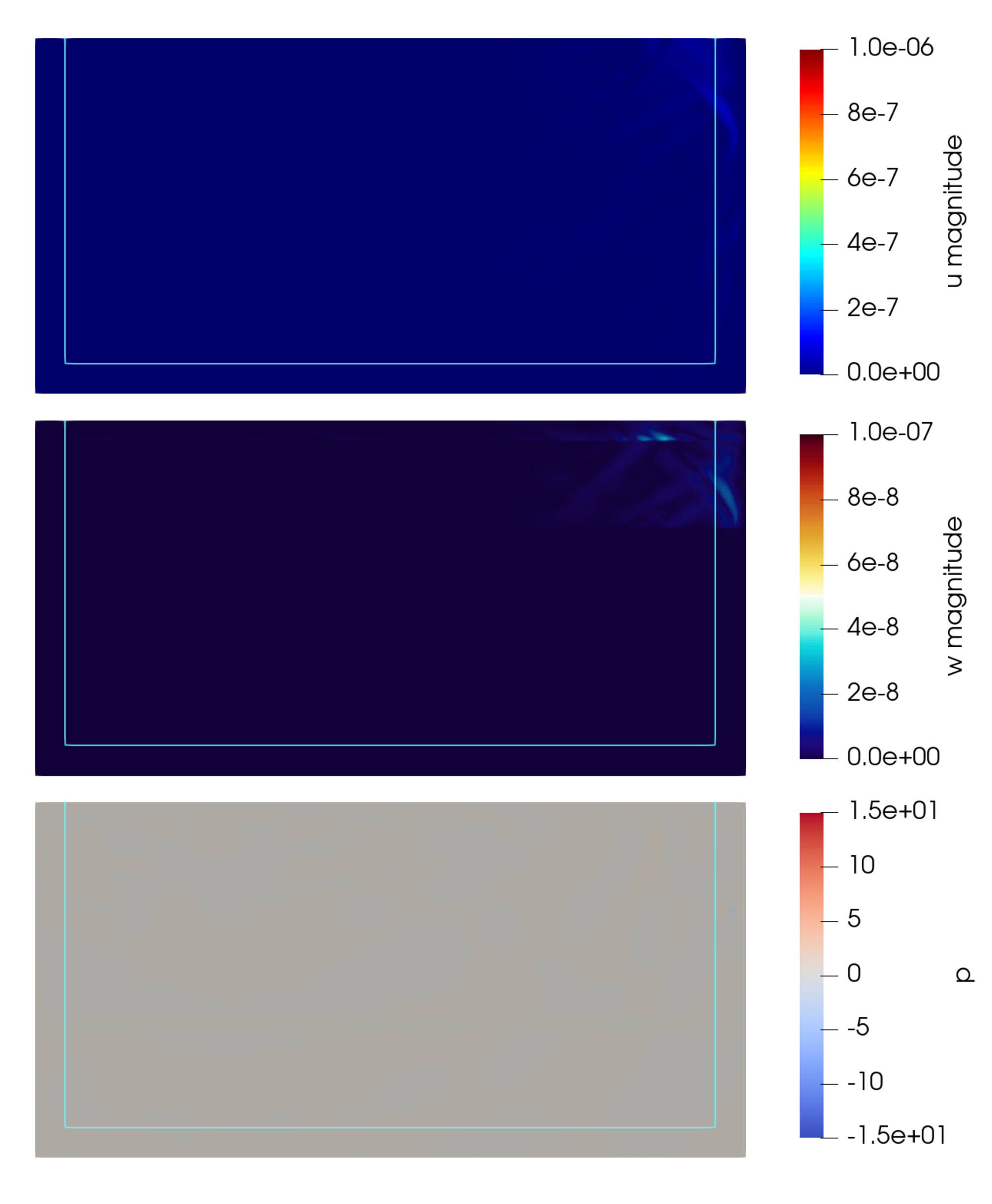} \label{fig: Exp 2 screenshots (c)}}

    \caption{Screenshots of the hybrid PML simulations at different time steps. The first row of the figure shows the solid displacement ($\uu$), the second the relative fluid displacement ($\ww$), and the third row the fluid pressure ($p$). No spurious reflections are observed at the interface between $\OmegaRD$ and $\OmegaPML$.}
    \label{fig: Exp 2 screenshots}
\end{figure}

Figure \ref{fig: heterogeneous energy a} shows the poroelastic energies for extended, paraxial, and PML simulations in the horizontally-layered domain. Similarly, the results indicate good agreement between the PML and reference solutions within the first 2 seconds of simulation, although some small differences are observed. The energy decay rate of the PML simulations is greater than that of the paraxial case, which decays slowly but consistently. In this experiment, it is not evident from the plots that outgoing waves are reaching the boundaries of $\OmegaRD$, as the plateaus observed in Figure \ref{fig: exp 1 energy} are absent. This is due to reflections at the interfaces between layers (see Figures \ref{fig: heterogeneous experiment} and \ref{fig: Exp 2 screenshots}).

\begin{figure}[h!]
    \centering
    \subfloat[Point A]{\includegraphics[width=0.49\textwidth]{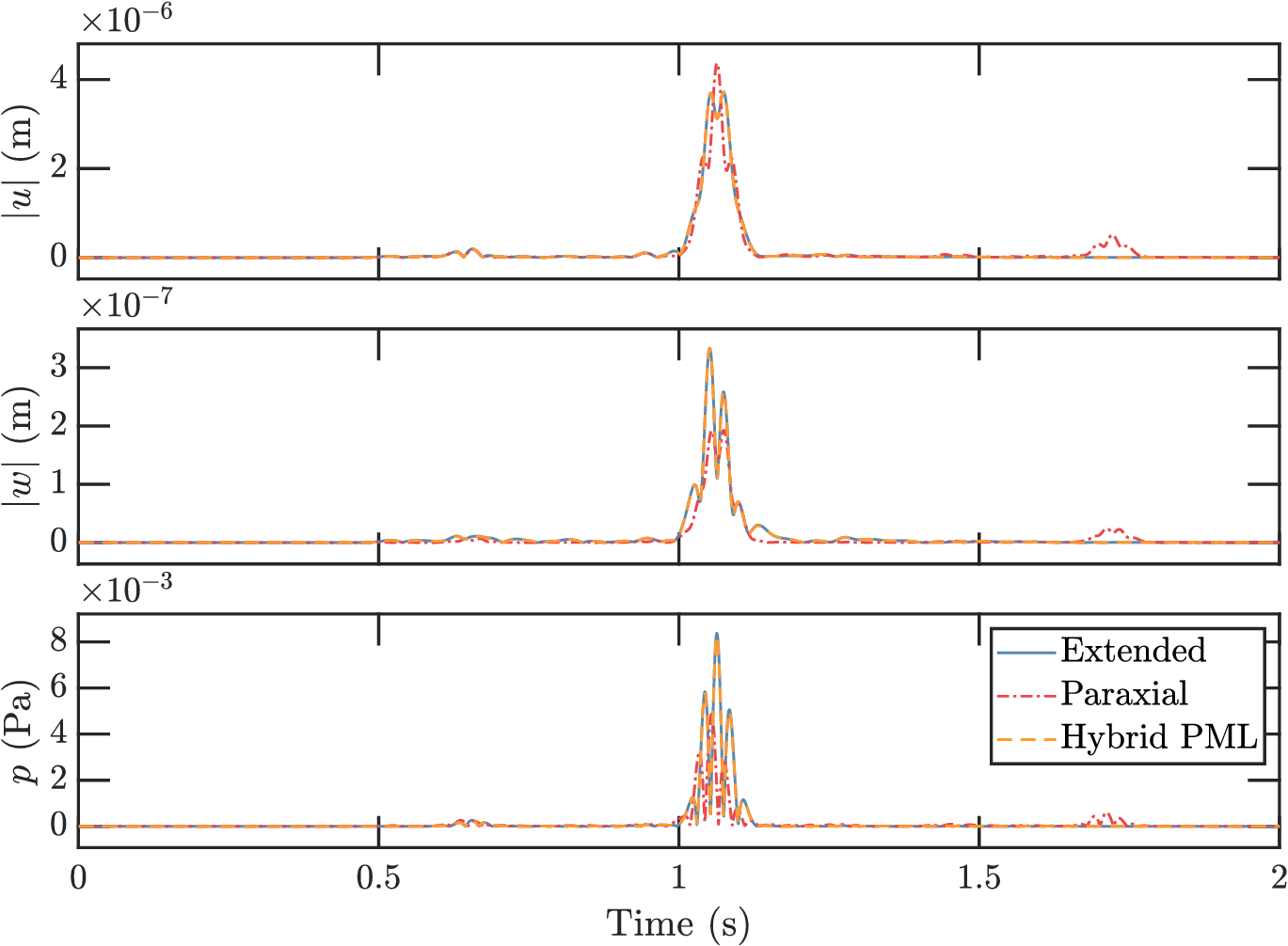}}\hfil
    \subfloat[Point B]{\includegraphics[width=0.49\textwidth]{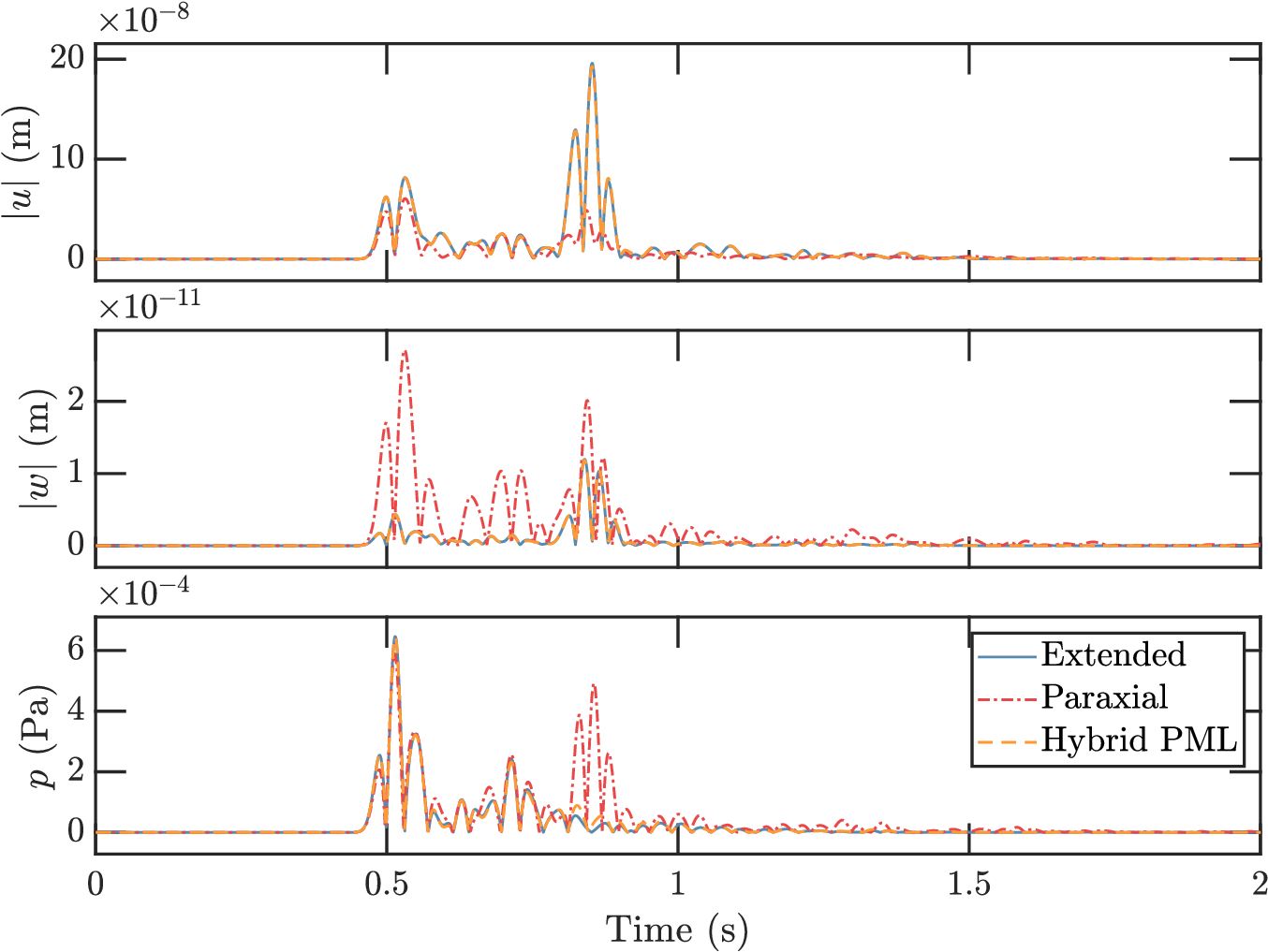}}

    \caption{Traces of $\uu$, $\ww$, and $p$ for the third experiment at the locations highlighted in Figure \ref{fig: heterogeneous experiment}. Results obtained using the hybrid PML formulation show good agreement with the extended domain solution.}
    \label{fig: traces 2nd experiment}
\end{figure}

The energies obtained using the hybrid and fully-mixed PML formulations produced almost the same results, illustrating that both methods provide similar solutions (see Figure \ref{fig: exp 1 energy a}). Additionally, the number of DOFs in the hybrid case is approximately 1.8 times fewer than the fully-mixed problem (as shown in Table \ref{tab3: degrees of freedom}). Therefore, the hybrid formulation is less computationally expensive than the fully-mixed form while providing equivalent solutions. Regarding Figure \ref{fig: heterogeneous energy b}, in terms of decay rate, during the first 2.5 seconds of the simulation, the results obtained using PML and M-PML in the hybrid formulation are similar, and slightly greater errors are observed with M-PML between 1.5 and 2.5s. The worst performance of M-PML in this time window is because the stretching functions are not in perfectly matched at $\GammaI$, as mentioned in previous paragraphs, which generates spurious reflections \cite{martin_high-order_2010}. However, during the last second, the energy of the uniaxial case decays more slowly than the M-PML case, which keeps a constant rate.

\begin{figure}[h!]
    \centering
    \subfloat[Point A]{\includegraphics[width=0.49\textwidth]{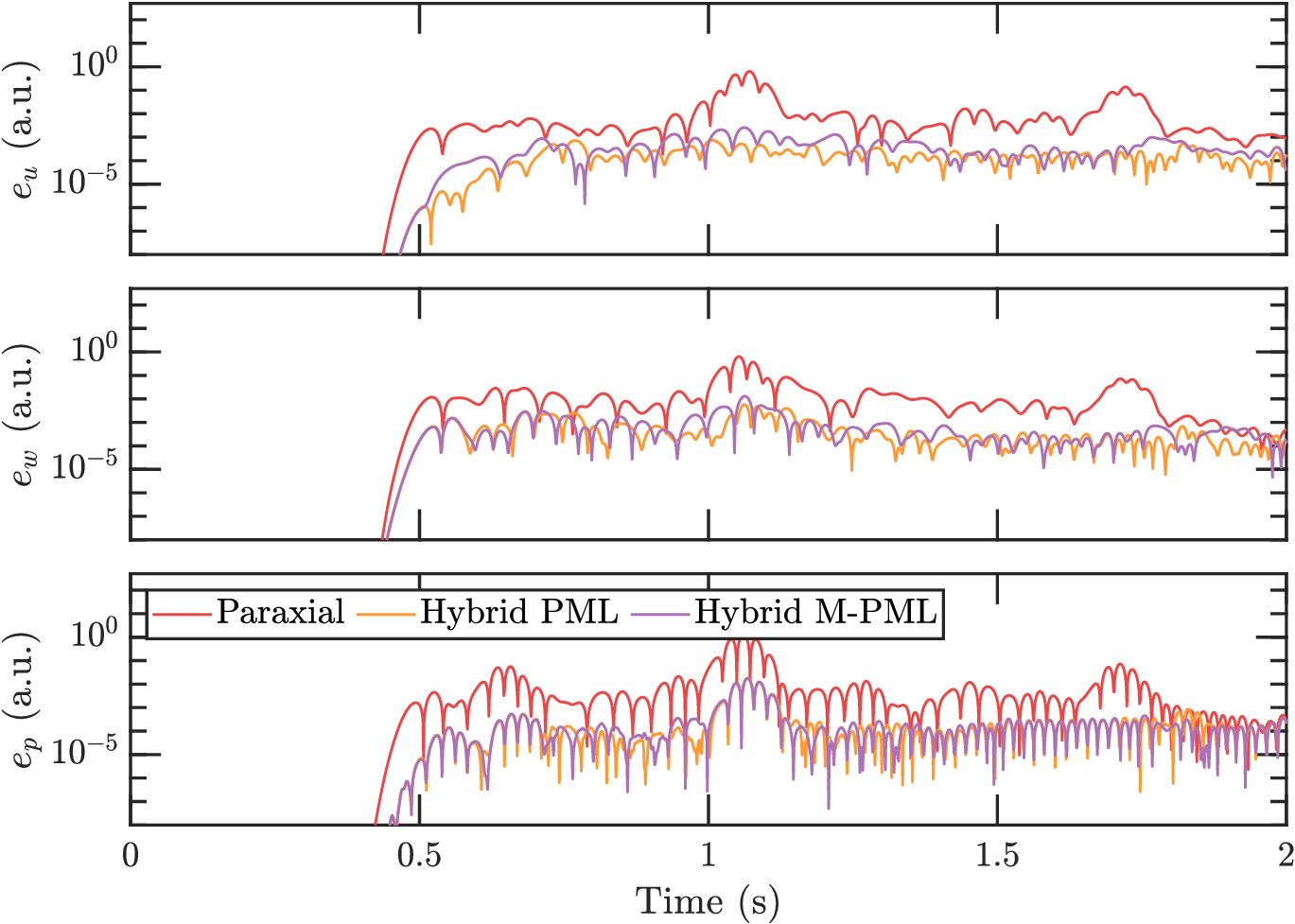}} \hfil
    \subfloat[Point B]{\includegraphics[width=0.49\textwidth]{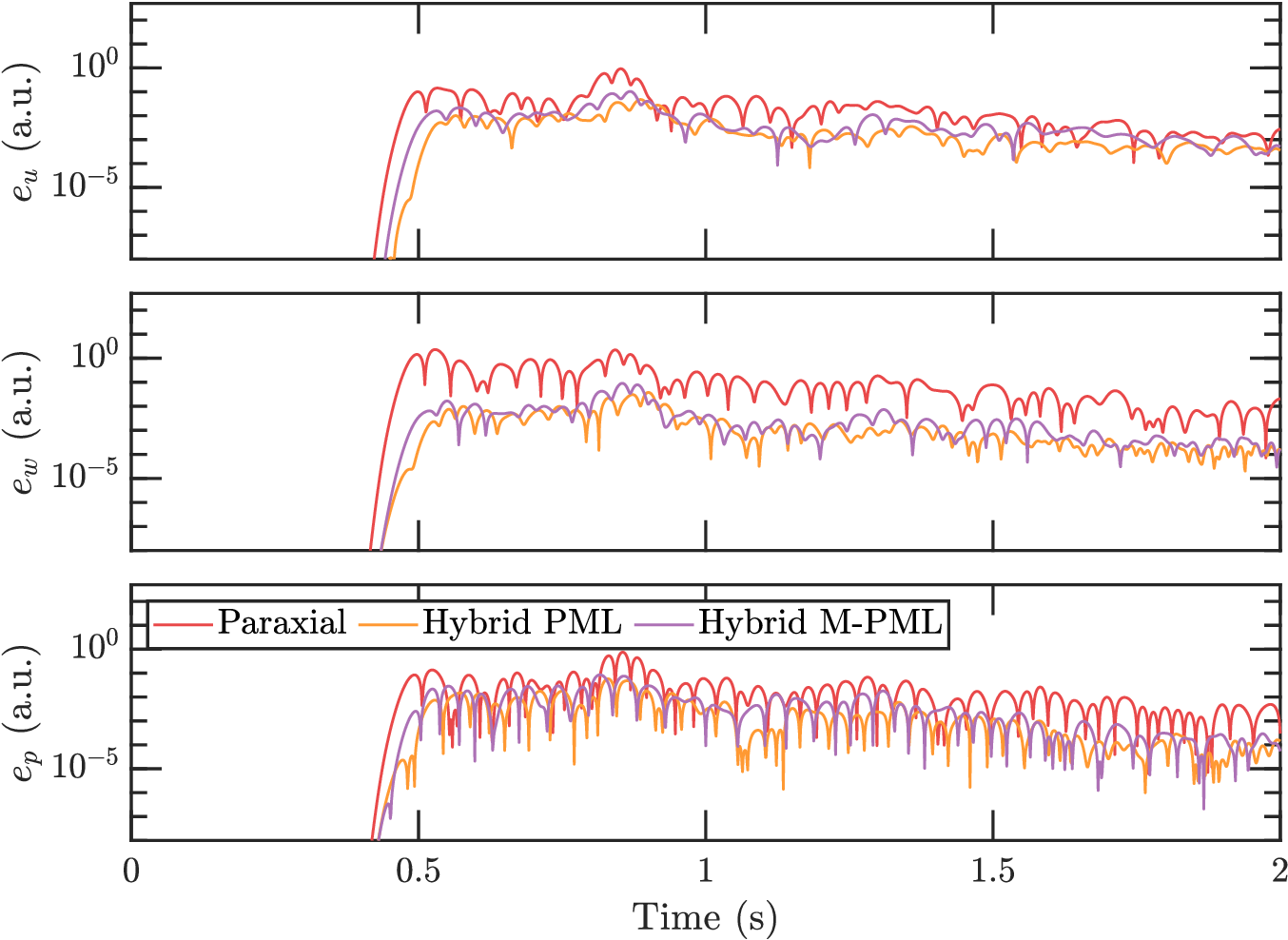}}

    \caption{Errors of the traces of $\uu$, $\ww$, and $p$ estimated using \eqref{eq: trace error} for the Experiment 2 at the locations highlighted in Figure \ref{fig: heterogeneous experiment}. The vertical axis is in logarithmic scale to improve differences visualization. At the beginning of the simulation, errors obtained with both methods are close to zero and therefore are out from the vertical axis range. Results obtained with the hybrid PML formulation do not show observable differences compared to the reference solution.}
    \label{fig: error of traces 2nd experiment}
\end{figure}

Figures \ref{fig: traces 2nd experiment} and \ref{fig: error of traces 2nd experiment} present traces of the solutions and errors calculated using \eqref{eq: trace error} for two locations shown in Figure \ref{fig: heterogeneous experiment}. The results show a good match between the hybrid PML and extended domain simulations at both locations, in contrast to the paraxial case where reflections are observed. The error analysis confirms the marginally superior performance of the hybrid PML method with respect to the M-PML case, but shows a considerable improvement compared to the paraxial case. The difference between PML and M-PML results is due to imperfect matching of stretching functions in the last case. Nevertheless, the results obtained with M-PML are still superior to those obtained with paraxial boundary conditions and more stable over time than the hybrid PML simulation (see Figure \ref{fig: heterogeneous energy}).

Finally, Figure \ref{fig: Exp 2 screenshots} presents screenshots that depict the propagation of waves at different times. These snapshots clearly show the transitions between different layers and the behavior of waves within each medium. For example, the layer with lower permeabilities (as shown in Table \ref{tab1: physical parameters} and Figure \ref{fig: heterogeneous experiment}) exhibited smaller relative fluid displacements ($\ww$) due to the small hydraulic conductivity. Additionally, media filled with air had lower pressures. Despite the complex behavior of waves in different media, the PML effectively absorbed and attenuated the waves during the simulation.

\subsection{Experiment 3: layered domain with outcropping}

\begin{figure}[h!]
    \centering
    \subfloat[]{\includegraphics[width=0.49\textwidth]{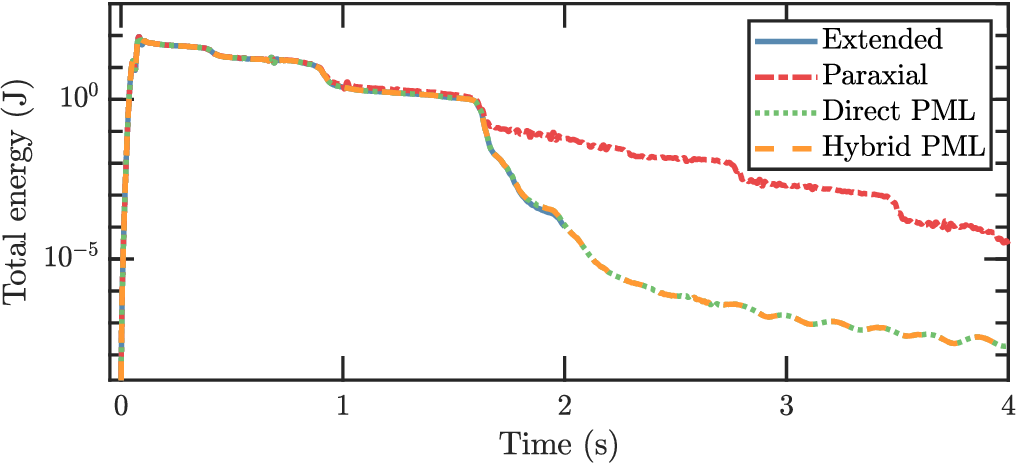} \label{fig: complex heterogeneous energy a}} \hfill
    \subfloat[]{\includegraphics[width=0.49\textwidth]{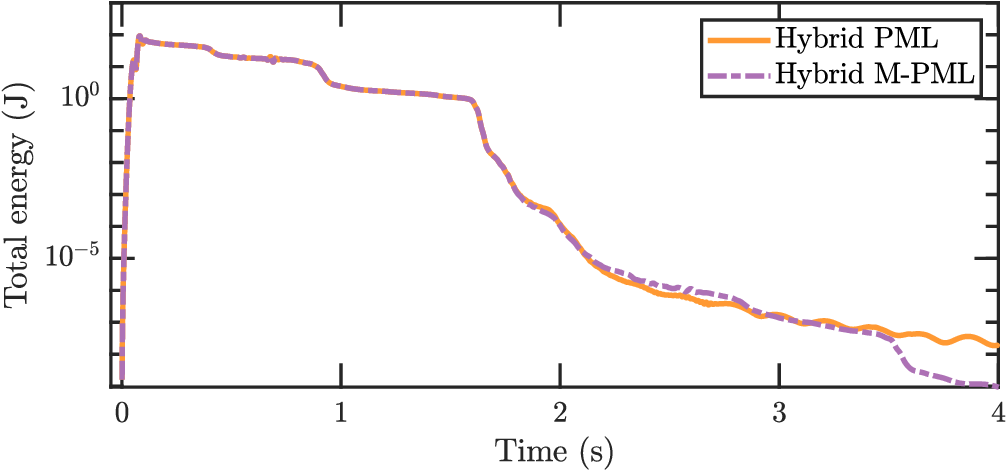} \label{fig: complex heterogeneous energy b}}
    \caption{Energy estimated on $\OmegaRD$ using  \eqref{eq: energy} for the Experiment 3. Results show the extended, paraxial, fully-mixed, and hybrid PML cases. Both fully-mixed and hybrid formulations give identical results.}
    \label{fig: complex heterogeneous energy}
\end{figure}

As observed in the previous experiments, both hybrid and fully-mixed PML formulations produced superior results compared to the paraxial case. The energy plots for this experiment also shows an excellent agreement between the PML and reference solutions within the initial 2 seconds of simulation, with negligible differences, as illustrated in Figure \ref{fig: complex heterogeneous energy a}. Although the decay of energy was consistent over time in all methods, PML exhibited faster energy decay compared to the paraxial case. The energy plots indicate that the solutions obtained using the hybrid and fully-mixed PML formulations are almost indistinguishable, indicating that both methods produce nearly identical results. This observation is consistent with the results of the previous experiments. Additionally, the hybrid case had a reduction in the number of DOFs by almost 1.8 times compared to the fully-mixed case (as shown in Table \ref{tab3: degrees of freedom}).

\begin{figure}[h!]
    \centering
    \subfloat[Location A]{\includegraphics[width=0.49\textwidth]{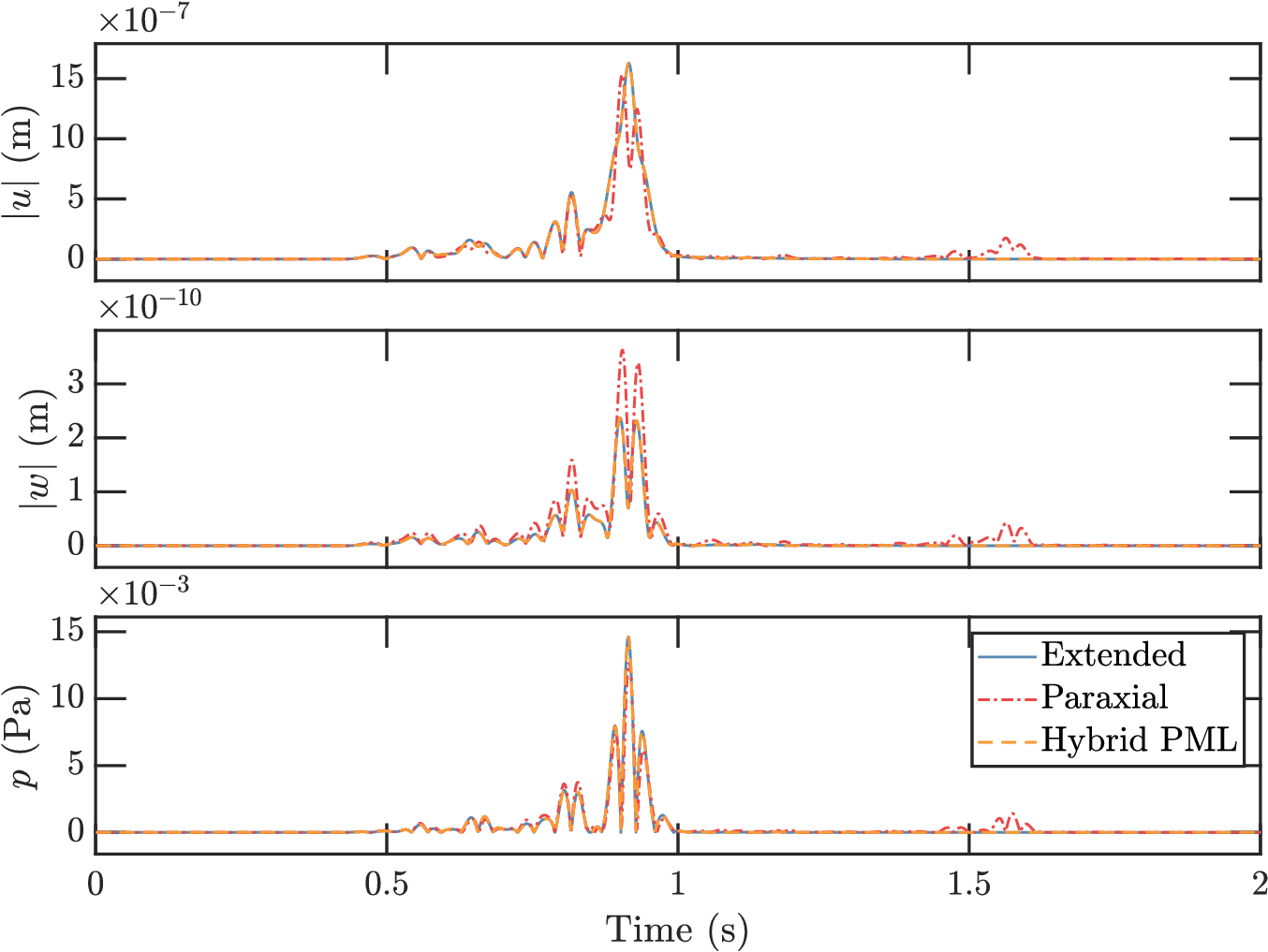}}\hfil
    \subfloat[Location B]{\includegraphics[width=0.49\textwidth]{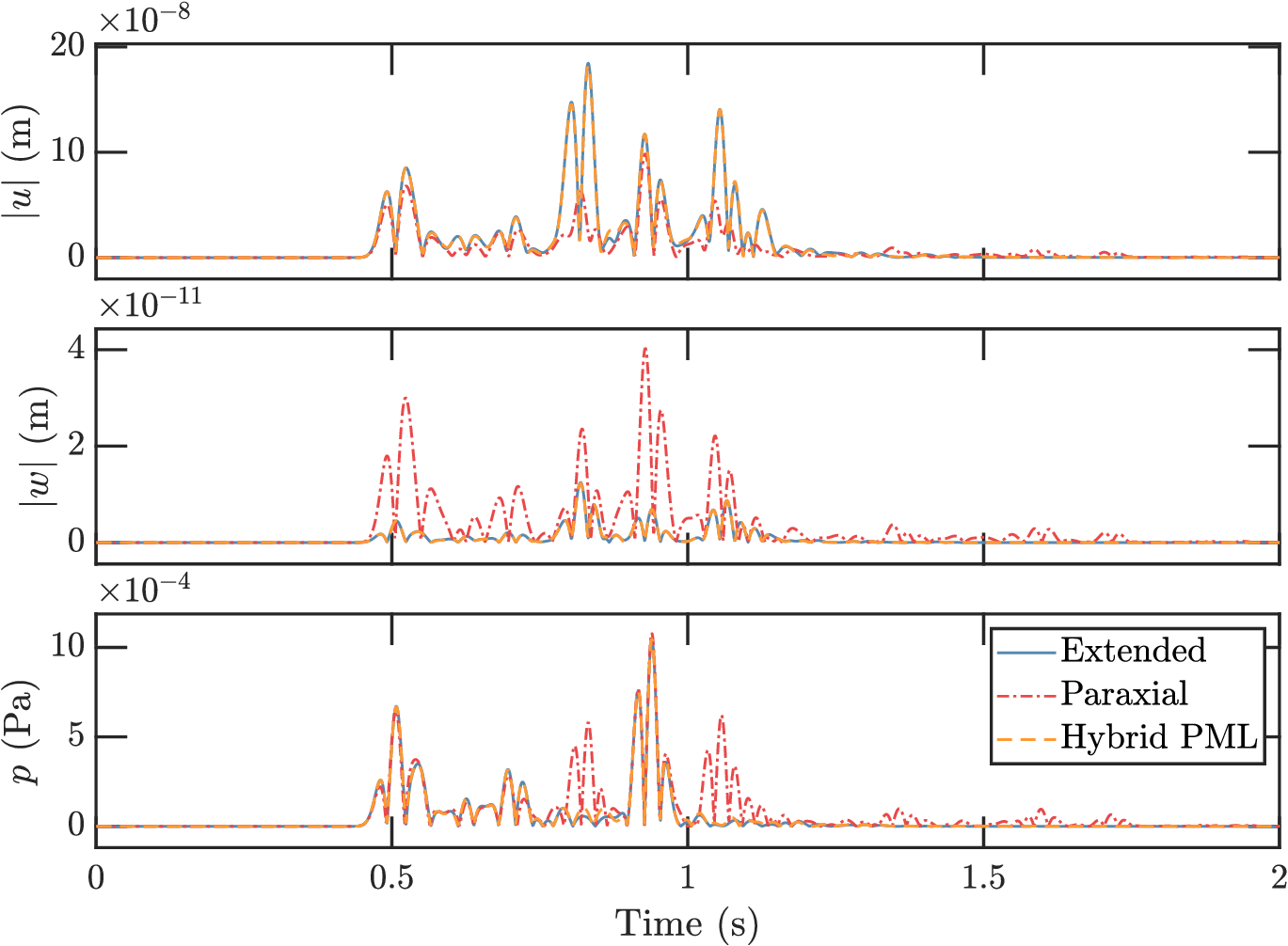}}

    \caption{Traces of $\uu$, $\ww$, and $p$ for the third experiment at the locations highlighted in Figure \ref{fig: complex heterogeneous experiment}. Results obtained using the hybrid PML formulation show good agreement with the reference solution.}
    \label{fig: traces 3rd experiment}
\end{figure}

Figure \ref{fig: complex heterogeneous energy b} compares the energies obtained using the hybrid PML and M-PML formulations. The results show that the two methods provide almost identical solutions up to seconds of the simulation, while some differences can be observed between 2 and 3 seconds, where PML slightly surpasses the multiaxial solution.. Interestingly, after 3 seconds of the simulation, the energy decay rate using M-PML is better with respect to the uniaxial case, highlighting the improved time-stability of M-PML stretching functions.

\begin{figure}[h!]
    \centering
    \subfloat[Point A]{\includegraphics[width=0.49\textwidth]{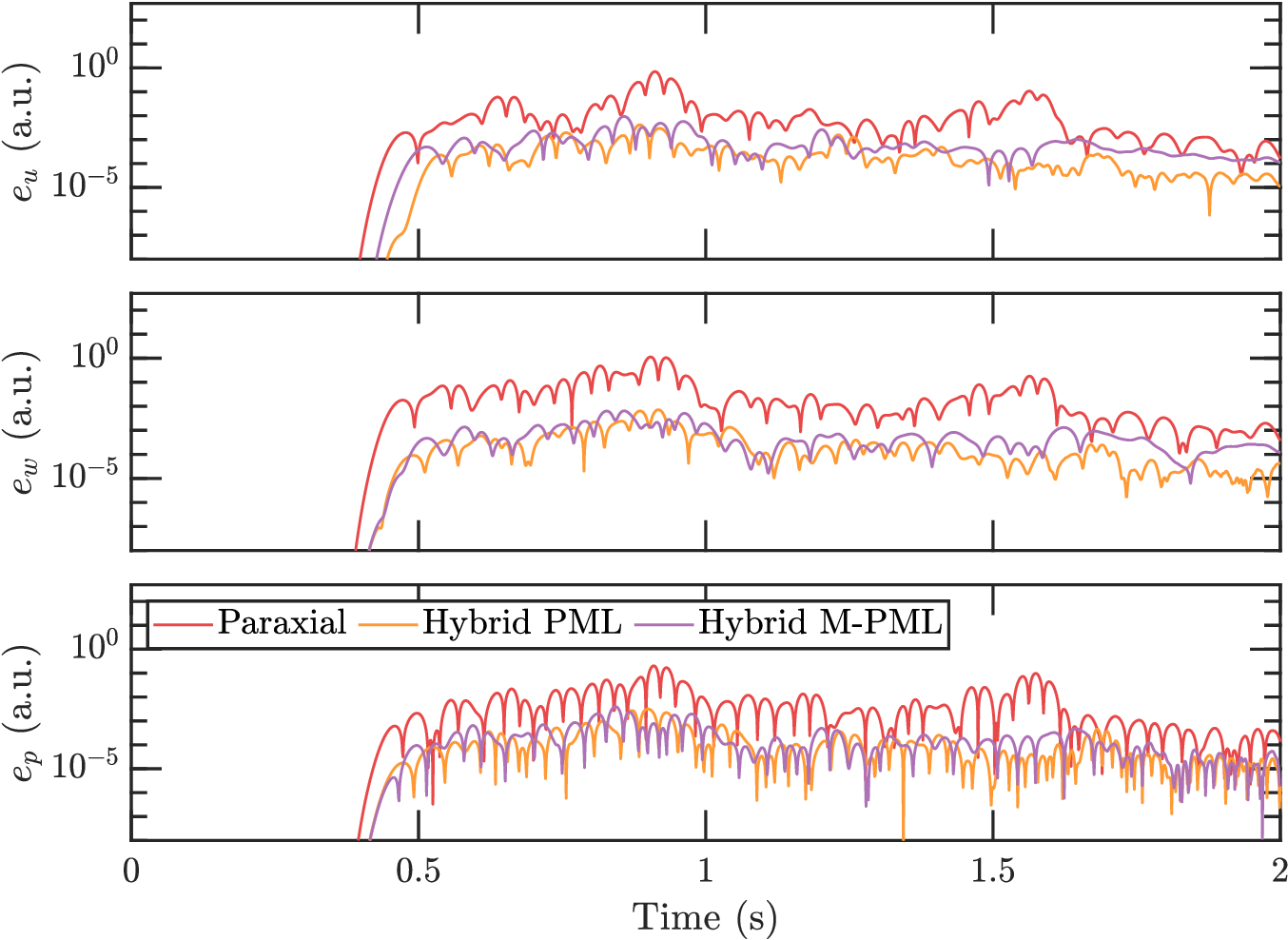}} \hfil
    \subfloat[Point B]{\includegraphics[width=0.49\textwidth]{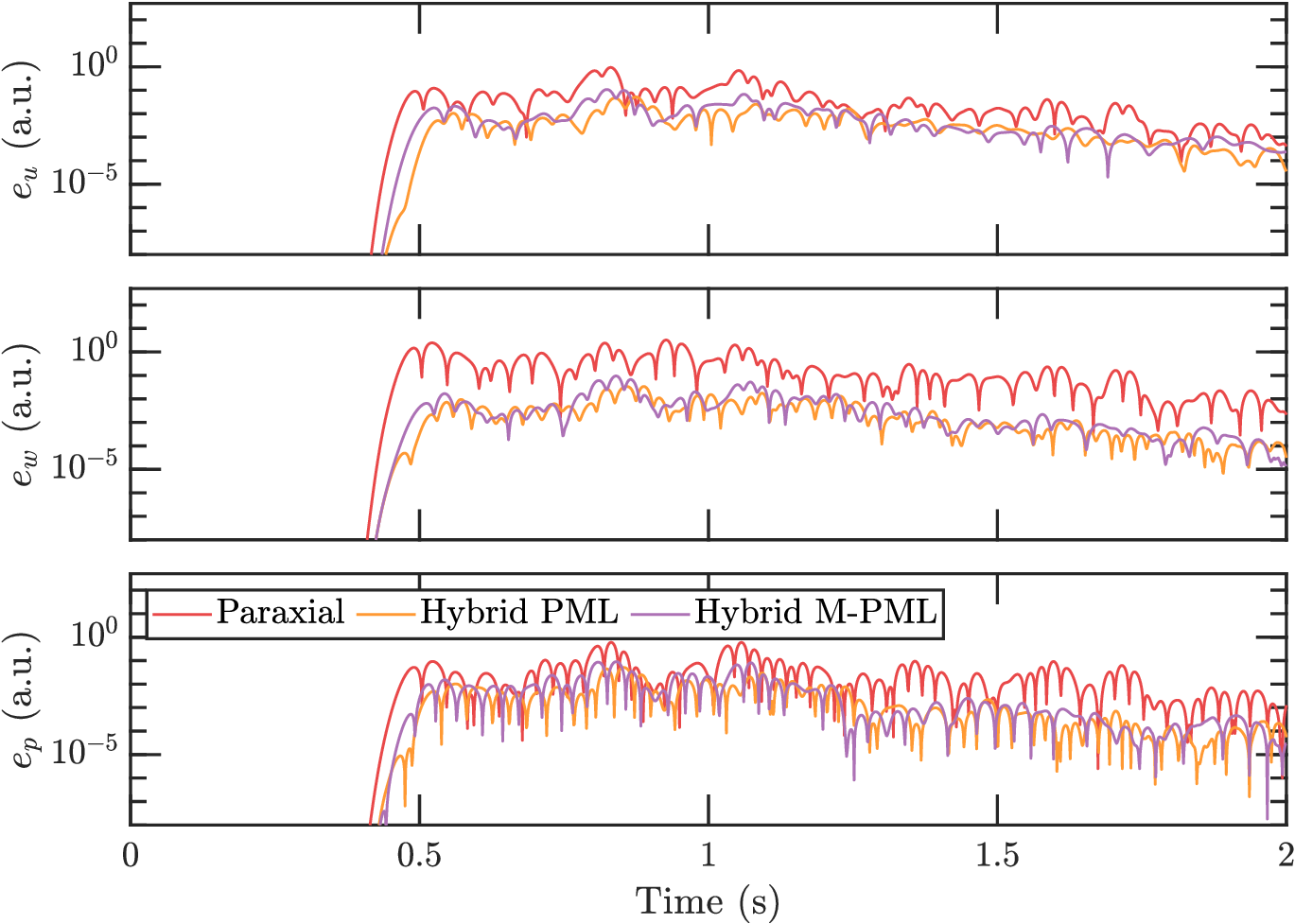}}

    \caption{Errors of the traces of $\uu$, $\ww$, and $p$ estimated using \eqref{eq: trace error} for the Experiment 3 at the locations highlighted in Figure \ref{fig: complex heterogeneous experiment}. Results obtained with the hybrid PML formulation do not show differences compared to the reference solution at Point A, while the paraxial approximation also performs reasonably well at this location.However, the advantages of the proposed method are much more evident at Point B. }
    \label{fig: error of traces 3rd experiment}
    
\end{figure}
\begin{figure}[h!]
    \centering
    \subfloat[$t=0.3$ s]{\includegraphics[scale=0.15]{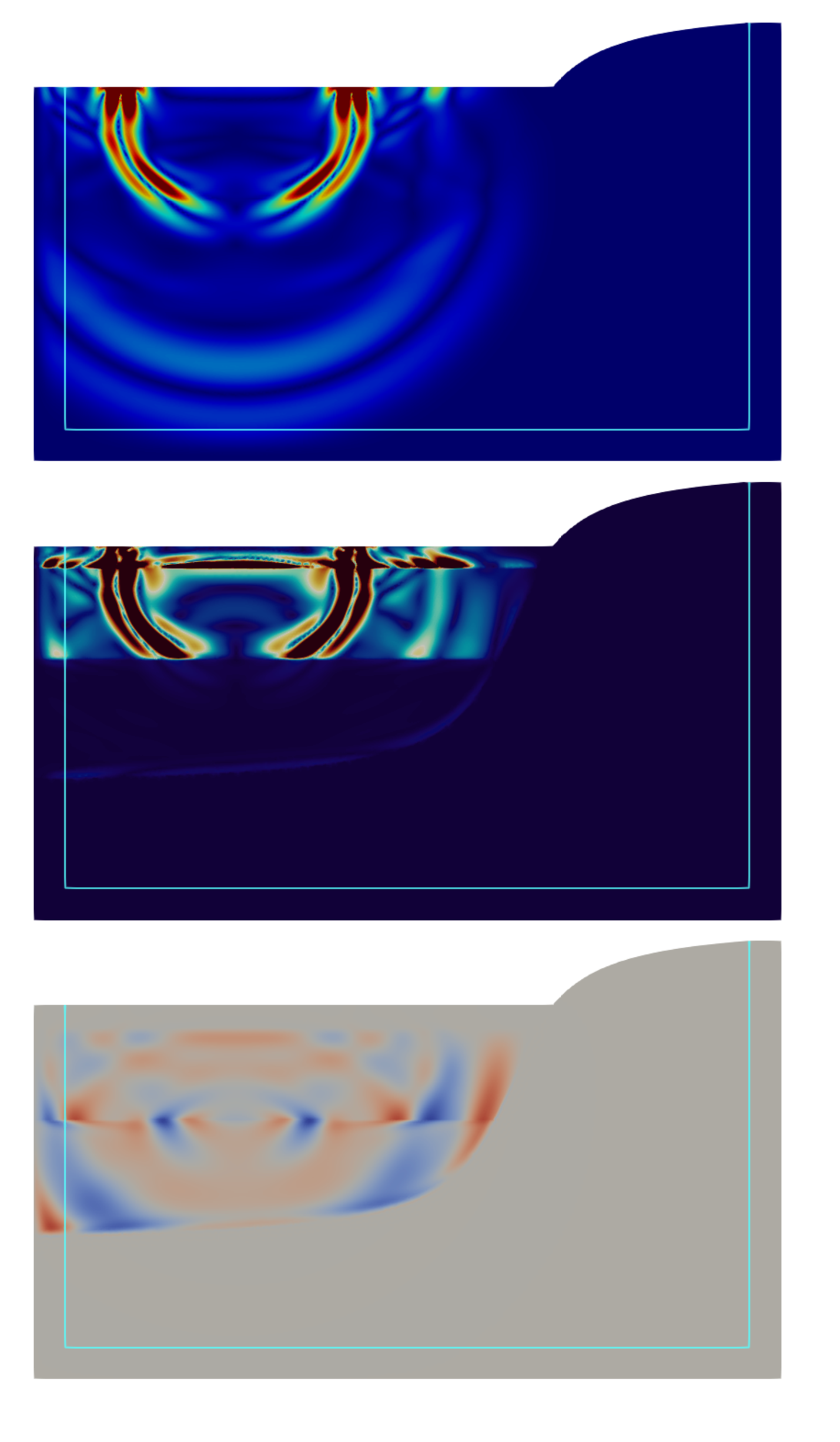} \label{fig: screenshots (a)}}
    \subfloat[$t=0.6$ s]{\includegraphics[scale=0.15]{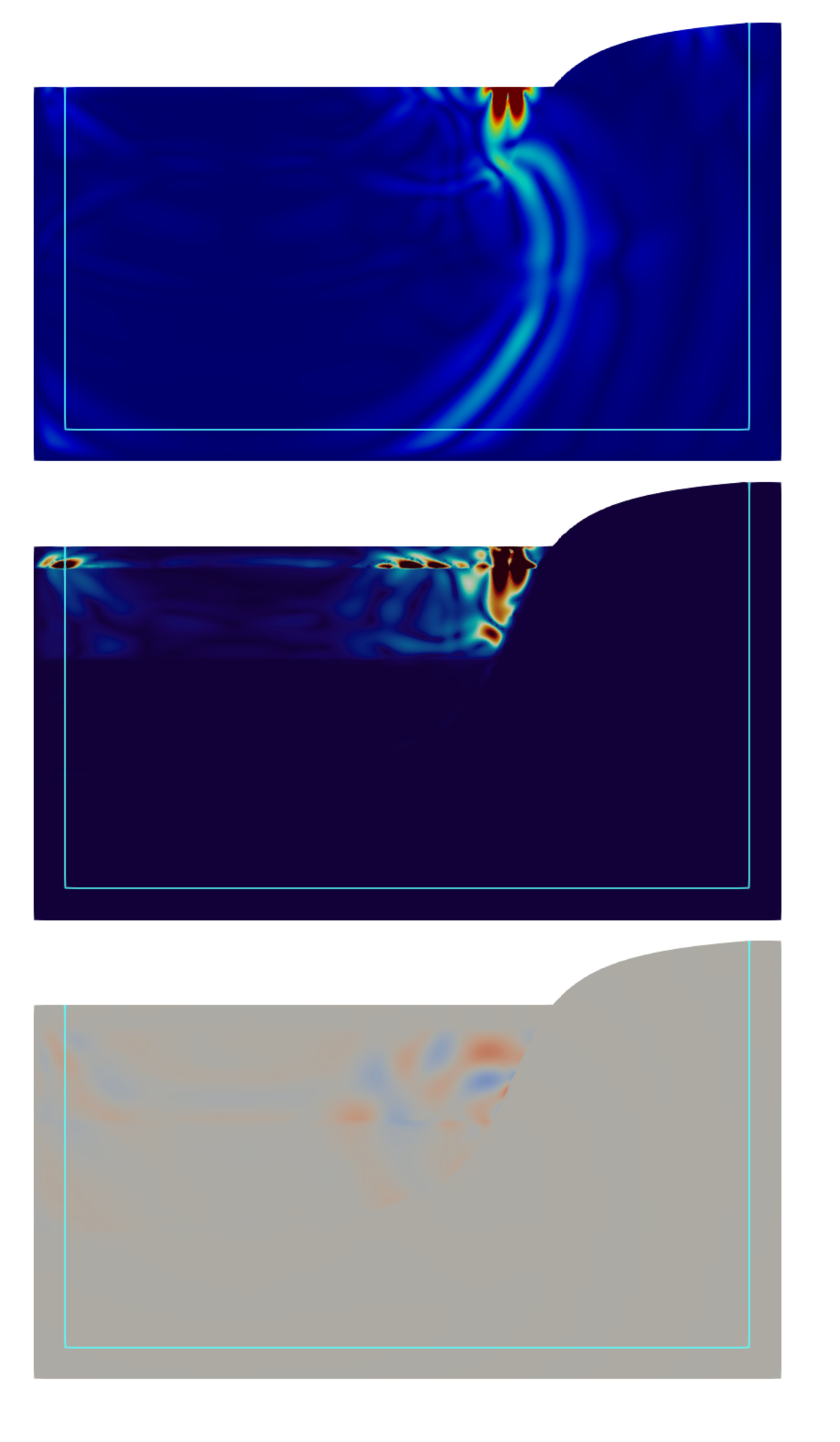} \label{fig: screenshots (b)}}
    \subfloat[$t=1.2$ s]{\includegraphics[scale=0.15]{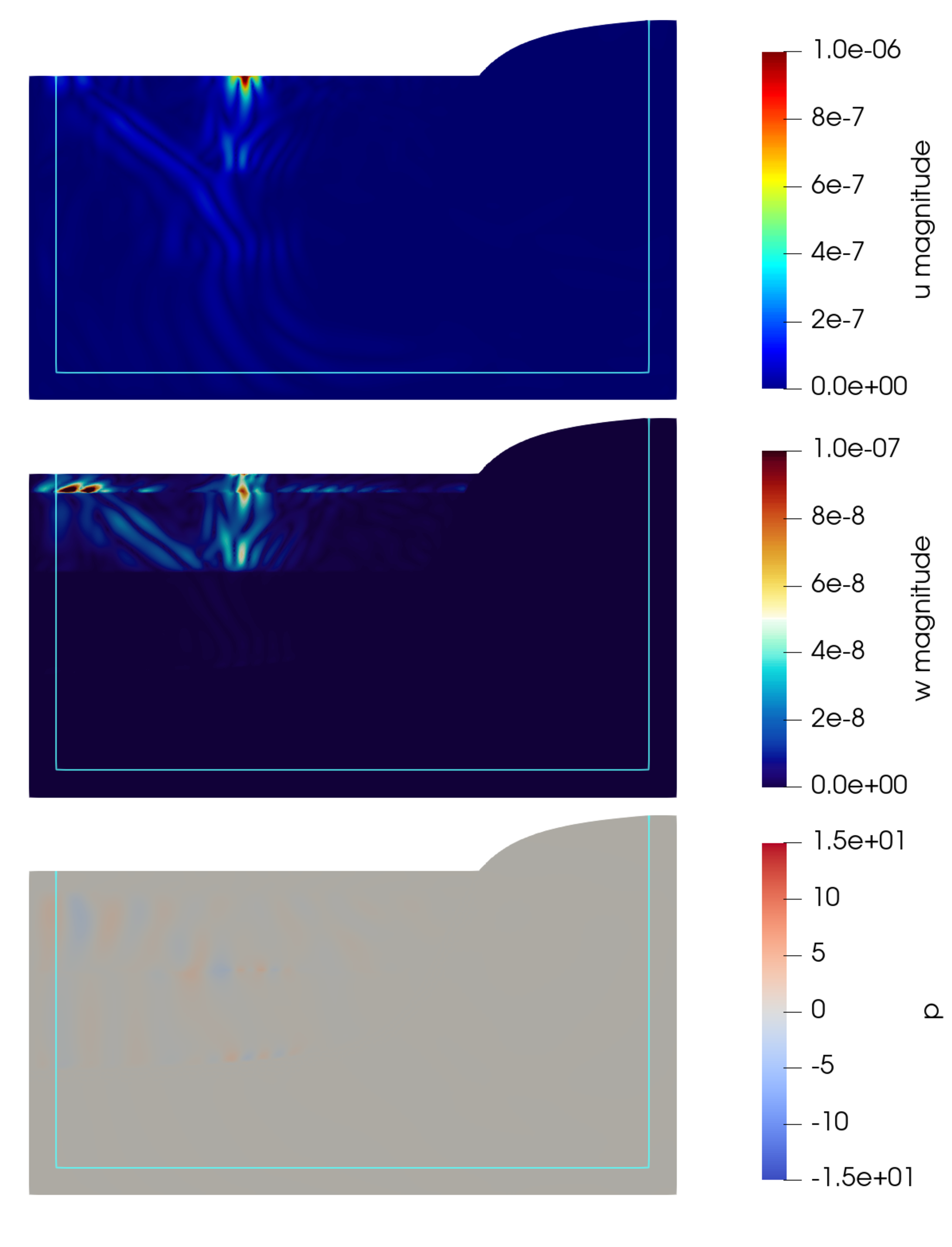} \label{fig: screenshots (c)}}

    \caption{Screenshots of the hybrid PML simulations at different time steps. The first row shows the solid displacement ($\uu$), the second the relative fluid displacement ($\ww$), and the third row the fluid pressure ($p$). No spurious reflections are observed at the interface between $\OmegaRD$ and $\OmegaPML$.}
    \label{fig: Exp 3 screenshots}
\end{figure}

Traces of the solutions and their corresponding errors, calculated using \eqref{eq: trace error}, for the two locations illustrated in Figure \ref{fig: complex heterogeneous experiment}, are presented in Figures \ref{fig: traces 3rd experiment} and \ref{fig: error of traces 3rd experiment}, respectively. These figures confirm the excellent agreement between the hybrid PML and the extended domain simulation at both locations, with no observable reflections unlike the paraxial case. The error analysis shows a considerable improvement compared to the paraxial case and slightly superior performance of the hybrid PML method with respect to the M-PML case. However, it is worth noting that the results obtained with M-PML are still better than those obtained with paraxial boundary conditions and exhibit better stability over time than the hybrid PML simulation.

Figure \ref{fig: Exp 3 screenshots} shows snapshots of the solutions at different times. Due to the complex  interfaces between materials in the medium, the waves interact in intricate ways generating complex patterns of reflected waves at interfaces as can be seen in the figure. Despite this complexity, the PML effectively absorbed the waves at the boundary of $\OmegaRD$, and no unwanted reflections were observed.




\section{Conclusions}

We proposed fully-mixed and hybrid formulations of the PML method for the simulation of poroelastic waves in truncated domains. Compared to other methods, both introduce only three additional scalar unknowns, i.e., the components of the symmetric stress-history tensor, reducing the number of unknowns with respect to current split- and unsplit-field formulations. However, the hybrid formulation considerably reduced the degrees of freedom required to solve the propagation problem when compared to the fully-mixed and extended domain simulations, because new unknowns are only defined in the boundary layer. On average, the hybrid PML formulation reduced the number of DOFs by approximately 1.8 times compared to the fully-mixed form and by 18 times compared to the extended domain simulations. However, in comparison to the paraxial case, the hybrid PML formulation increased the number of DOFs by almost 41\%. However, although effective for some applications, paraxial boundary conditions are not ideal for absorbing surface waves. Therefore, they are not recommended in problems where surface waves are predominant.

The proposed formulations of the PML method are prone to the same issues of other PML formulations, and they may suffer instability over time under certain circumstances. However, a significant advantage of the proposed methods is that they enable the scaling and attenuation functions to be redefined using stretching functions with superior absorbing properties, such as M-PML, without modifying the underlying partial differential equations. Moreover, the time-integration scheme used for the simulations did not consider numerical damping, making the problem conditions even more demanding compared to other methods \cite{He2019a,He2019}. Only small time instabilities were observed, and they were fixed using M-PML.

In terms of discretization, the element size was selected as large as possible to reduce the number of DOFs. Additionally, only P2-P2-P1 finite element triplets were utilized to represent the solutions to the problems. This approach considerably reduced the computational effort required to solve the problem, unlike in previous studies. As a result, the hybrid and fully-mixed PML formulations proposed in this article have demonstrated robustness for demanding discretization and physical media conditions.


The following steps regarding this investigation are: (1) to extend fully-mixed and hybrid formulations to the 3D case to simulate more realistic scenarios for exmaple for complex seismic geophysical applications. (2) To consider variable and discontinuous porosities in space which would introduce discontinuities in the relative fluid displacement. (3) To apply the 2D and 3D formulations to solve inverse problems in porous media. And (4), to simulate and understand the propagation of poroelastic waves in the human body, which is particularly interesting in the fields of Magnetic Resonance Imaging and Ultrasound and is closely related to the elastography problem \cite{doyley_model, Kolipaka2009}.

\section{Acknowledgments}
HM and JM acknowledge the financial support given by ANID through the projects ANID-FONDECYT Postdoctorado \#3220266 and ANID-FONDECYT Regular \#1230864, respectively. ES was partially funded by a grant from the Research Center for Integrated Disaster Risk Management CIGIDEN Project 1522A0005 FONDAP 2022.


\bibliographystyle{elsarticle-harv}


\end{document}